\documentclass[english, 11pt]{amsart}

\usepackage{tikz}
\usepackage{aliascnt}
\usepackage{mathrsfs}
\usepackage[all,poly,knot]{xy}
\usepackage{hyperref}
\usepackage{comment}  
\usepackage{csquotes}   
\usepackage{amssymb,amsbsy,amsmath,amsfonts,amssymb,amscd,
	graphics,color,footmisc,fancyhdr,multicol,fancybox,
	graphicx,mathrsfs,rotating,ifthen,wasysym}

\usepackage{lipsum}

\usepackage{times}
\usepackage{mathtools}
\usepackage{bm}
\usepackage{pdfpages}
\usepackage{multirow}
\usepackage{nccrules}
\usepackage{textcomp}
\usepackage{comment}
\usepackage{framed}
\usepackage[all]{xy}
\usepackage{graphicx}
\usepackage{pgf,tikz}
\usepackage{mathrsfs}
\usetikzlibrary{arrows}
\usepackage{lipsum}
\usepackage{mathtools}

\DeclareMathOperator{\dist}{dist}

\DeclareMathOperator{\Id}{Id}
\DeclareMathOperator{\FS}{FS}
\DeclareMathOperator{\End}{End}

\newcommand{\C}{\mathbb{C}}

\newcommand{\D}{\mathbb{D}}

\def\log{\mathrm{log}\,}

\theoremstyle{plain}

\newtheorem{thm}{Theorem}[section]

\newtheorem{lem}[thm]{{Lemma}}

\newtheorem{pro}[thm]{Proposition}

\newtheorem*{claim}{Claim}

\theoremstyle{remark}
\newtheorem{rmk}[thm]{Remark}

\numberwithin{equation}{section}

\usepackage[top=2cm, bottom=2cm, left=2.5cm, right=2.5cm]{geometry}

\usepackage[hyperpageref]{backref}

\usepackage{xcolor}
\hypersetup{
	colorlinks,
	linkcolor={red!50!black},
	citecolor={blue!62!black},
	urlcolor={blue!80!black}
}
\theoremstyle{plain}
\newcommand{\thistheoremname}{}
\newtheorem*{genericthm*}{\thistheoremname}
\newenvironment{namedthm*}[1]{\renewcommand{\thistheoremname}{#1}%
	\begin{genericthm*}}
	{\end{genericthm*}}

\newtheoremstyle{named}{}{}{\itshape}{}{\bfseries}{.}{.5em}{\thmnote{#3's }#1}
\theoremstyle{named}

\makeatletter
\newcommand\thankssymb[1]{\textsuperscript{\@fnsymbol{#1}}}
\makeatother
\begin{document} 
	\title[On the value distribution theory of $J$-holomorphic Brody curves]{On the value distribution theory\\
		of $J$-holomorphic Brody curves}

	
\author{Dinh Tuan Huynh}

\address{Department of Mathematics, University of Education, Hue University, 34 Le Loi St., Hue City, Vietnam}
\email{dinhtuanhuynh@hueuni.edu.vn}
	
\author{Yunling Chen}
\address{Academy of Mathematics and Systems Sciences, Chinese Academy of Sciences, Beijing 100190, China}
\email{chenyl25@amss.ac.cn}

\begin{abstract}
Let $J$ be a smooth almost complex structure on $\mathbb{CP}^2$ tamed by the Fubini--Study form. For a $J$-holomorphic Brody curve $f:\mathbb{C}\to (\mathbb{CP}^2,J)$ and a family of $q\ge2$ $J$-lines $\{L_i\}_{1\le i\le q}$ in general position, we prove the following Second Main Theorem type estimate
\[
(q-1)T_f(r)\le \sum_{i=1}^q N_f(r,L_i)+o(r^2).
\]
This extends a result of Da Costa-Duval \cite{Dacosta-Duval2013} to almost complex setting. As an application, we show that if a Brody $J$-holomorphic curve \(f\) omits two $J$-lines, then its area has linear energy growth
\[
\int_{\D_r}f^*\omega_{\FS}=O(r),
\]
which answers a question posted by Gaussier--Tsukamoto \cite{Gaussier-Tsukamoto}.
\end{abstract}

\maketitle

\section{Introduction}
Let \(X\) be a compact complex manifold with a hermitian metric \(\|\cdot\|\). A \emph{Brody curve} in \(X\) is an entire curve \(f\colon\mathbb{C}\to X\) whose derivative \(\|\mathrm{d}f\|\) is bounded. Equivalently, \(f\) is a Brody curve if and only if the family of its translations \(\{z\mapsto f(z+a):z\in\mathbb{C}\}\) is normal. These curves arise as limits of sequences of holomorphic maps from the unit disc whose derivatives at the origin blow up, after reparametrization and passage to a subsequence \cite{Brody78}. The absence of Brody curves is an important characterization of Kobayashi hyperbolicity for compact complex manifolds.

The set of all Brody curves is itself a natural infinite-dimensional dynamical system. Let \(X\) be compact, and let \(\mathcal{M}(X)\) be the space of Brody curves \(f:\mathbb{C}\to X\) with a fixed derivative bound, endowed with the compact-open topology. Then \(\mathbb{C}\) acts on \(\mathcal{M}(X)\) by translations
\[
(a\cdot f)(z)=f(z+a).
\]
Gromov introduced mean dimension precisely to measure the size of dynamical systems of this type, and he initiated the study of spaces of holomorphic maps from noncompact sources \cite{Gromov99}.

When \(X=\mathbb{CP}^n\), the standard complex projective space with the Fubini--Study metric \(\omega_{\mathrm{FS}}\), the Brody property for an entire curve \(f\colon\mathbb{C}\to\mathbb{CP}^n\) means that its spherical derivative
\[
\|\mathrm{d}f\|=\frac{\sqrt{\sum_{i<j}|f_i'f_j-f_j'f_i|^2}}{\|f\|^2}
\]
is bounded, where \([f_0:f_1:\dots:f_n]\) is a homogeneous presentation of \(f\) and \(\|f\|^2=\sum_{i=1}^{n}|f_i|^2\). The study of the mean dimension of \((\mathcal{M}(\mathbb{CP}^n):\mathbb{C})\) has led to a remarkably precise picture, thanks to enormous efforts by Tsukamoto and his collaborators \cite{Tsukamoto08}, \cite{Tsukamoto09a}, \cite{Tsukamoto09b}, \cite{Tsukamoto12}, \cite{Tsukamoto15}, \cite{Tsukamoto18}. These efforts culminate in the formula
\[
\dim\bigl(\mathcal{M}(\mathbb{CP}^n):\mathbb{C}\bigr)
=2(N+1)\rho(\mathbb{CP}^n),
\]
where
\[
\rho(\mathbb{CP}^n)
:=
\sup_{f\in\mathcal{M}(\mathbb{CP}^n)}\rho(f)
=
\sup_{f\in\mathcal{M}(\mathbb{CP}^n)}
\biggl(\limsup_{r\to\infty}\frac{1}{\pi r^2}\int_{|z|<r}f^*\omega_{\mathrm{FS}}
\biggr)
\]
is the maximal energy density of Brody curves in \(\mathbb{CP}^n\). This result shows explicitly that the large-scale energy of individual Brody curves controls the topological size of the whole dynamical system.

Quantitative growth estimates for the areas of Brody discs also play an important role in value distribution theory. Recall that the Ahlfors--Shimizu order function of a holomorphic curve \(f\colon\mathbb{C}\to\mathbb{CP}^n\) is defined by
\[
T_f(r)=\int_1^r\frac{\mathrm{d}t}{t}\int_{|z|<t}f^*\omega_{\mathrm{FS}}.
\]
Thus, if \(f\) is Brody, its order function grows at most like \(O(r^2)\). The first remarkable result on the growth rates of Brody curves is the classical theorem of Clunie and Hayman \cite{Clunie-Hayman66} (see also \cite{Pommerrenke70} for another proof). In modern language, it states that if a Brody curve \(f:\mathbb{C}\to\mathbb{CP}^1\) omits one point, then its order function has low growth that
\[
T_f(r)=O(r).
\]
Eremenko \cite{Eremenko-2010} generalized this phenomenon to higher-dimensional projective space: if a Brody curve \(f:\mathbb{C}\to\mathbb{CP}^n\) omits \(n\) hyperplanes in general position, then \(T_f(r)=O(r)\). Consequently, the mean dimension of the space of Brody curves into \(\mathbb{CP}^n\) that avoid \(n\) hyperplanes in general position is zero. Da Costa and Duval \cite{Dacosta-Duval2013} subsequently introduced a potential-theoretic method for Brody curves and established the following Second Main Theorem type estimate:
\[
(q-n+1)T_f(r)
\le \sum_{i=1}^q N_f(r,H_i)+o(r^2),
\]
for \(q\) hyperplanes \(\{H_i\}_{1\le i\le q}\) in general position. By allowing the error term to be \(o(r^2)\), this result improves significantly Cartan's Second Main Theorem, in which the coefficient bound \(q-n-1\) is required.

The concept of hyperbolicity has also been generalized to the almost-complex setting \cite{Boris-Marius99}, \cite{Kobayashi01}, \cite{Fathi-Adel09}. The hyperbolicity of certain configurations in \(\mathbb{CP}^2\) has been extended to the almost-complex case by \cite{Duval-Fourier04}, \cite{Saleur14}. Notably, on the quantitative side, a variant of Cartan's SMT of holomorphic functions in \(\mathbb{CP}^2\) has been developed in the setting of elliptic projective planes \cite{Duval-GAFA}.

Recently, Gaussier and Tsukamoto \cite{Gaussier-Tsukamoto} studied the energy of quasiconformal maps and pseudoholomorphic Brody curves in almost-complex setting. Here a pseudoholomorphic curve is taken with respect to a given almost-complex structure \(J\) on \(\mathbb{CP}^2\), and a \(J\)-line means a pseudoholomorphic line with respect to \(J\). Among other results, they proved that a Brody \(J\)-holomorphic curve in \(\mathbb{CP}^2\) omitting three \(J\)-lines in general position has zero energy density. Motivated by the classical linear growth  result in \cite{Eremenko-2010} for Brody curves, they asked \cite[Problem~0.4(2)]{Gaussier-Tsukamoto} whether omitting two \(J\)-lines already forces the sharp estimate
\[
\int_{|z|<r}f^*\omega_{\FS}=O(r).
\]

The first purpose of this paper is to establish an almost-complex analogue of the Da Costa--Duval Second Main Theorem Main Theorem \cite[Theorem 3]{Dacosta-Duval2013} in \(\mathbb{CP}^2\). Our first result states as follows.
\begin{thm}\label{thm:main}
	Let $J$ be a smooth almost complex structure on $\mathbb{CP}^2$ tamed by the Fubini--Study form $\omega_{\FS}$. Let $f:\mathbb{C}\longrightarrow (\mathbb{CP}^2,J)$be a Brody $J$-holomorphic curve. Let $L_1,\dots,L_q$, $q\ge2$, be $J$-lines in general position. If the image of $f$ is not contained in any $L_i$, then
	\begin{equation}\label{eq:main-smt}
	(q-1)T_f(r)\le \sum_{i=1}^q N_f(r,L_i)+o(r^2).
	\end{equation}
\end{thm}

Here and below, general position means that no three of the $J$-lines have a common point. It follows immediately from Theorem \ref{thm:main} that under the assumptions therein, if the Brody $J$-curve $f$ omits two $J$-lines in general position, then both the area \(a_f(r)\coloneq\int_{|z|<r}f^*\omega_{\FS}\) and the order function $T_f(r)$ have the growth rate as $o(r^2)$. This is sufficient to conclude that the energy density of $f$ is zero, but nevertheless not the optimal one suggested in \cite[Problem~0.4(2)]{Gaussier-Tsukamoto}. 

The second purpose of our paper is to prove the sharp quantitative result and thereby provide a positive answer to the above question of Gaussier--Tsukamoto.
\begin{thm}
\label{thm:gt-linear}
	Let $J$ be a smooth almost complex structure on $\mathbb{CP}^2$ tamed by $\omega_{\FS}$, and let $f:\mathbb{C}\longrightarrow (\mathbb{CP}^2,J)$	be a Brody $J$-holomorphic curve. If $f$ omits two distinct $J$-lines, then
	\begin{equation}\label{eq:gt-linear-characteristic}
	T_f(r)=O(r),
	\end{equation}
	and
	\begin{equation}\label{eq:gt-linear-area}
	a_f(r)=\int_{|z|<r}f^*\omega_{\FS}=O(r).
	\end{equation}
\end{thm}

We now briefly explain the key ideas in the proofs of our main results. To prove Theorem~\ref{thm:main}, for each \(J\)-line \(L_i\), let \(u_i\) denote its logarithmic potential. The almost-complex Poincar\'e--Lelong formula then takes the form
\[
\Delta u_i=f^*[L_i]+\operatorname{Div}F_i,
\]
where \(\|F_i\|_\infty\le C\) uniformly. Note that in the standard complex situation, the divergence term $\operatorname{Div}F_i$ is absent. In the almost complex case, we first solve away the bounded divergence error with zero boundary values (see Lemma \ref{lem:perturbed}). After redoing all the local computations and potential-theoretic results of \cite[Theorem 2]{Dacosta-Duval2013}, and taking into account all local error terms, we obtain
\[
M_r(u)=M_r\bigl(\max\{u_i,u_j\}\bigr)+O(r\log r),
\]
where
\[
M_r(u):=\frac1{2\pi}\int_0^{2\pi}u(re^{i\theta})\,\mathrm{d}\theta.
\]
The error term we receive at the end is $O(r\log r)$, instead of $O(r)$ in the usual complex setting. But this term is nevertheless $o(r^2)$, which is still enough for the proof of Theorem \ref{thm:main}. Indeed, after taking quadratic renormalization, the first-order perturbations disappear, and the zero-set argument of Da Costa--Duval, together with Grishin's lemma, yields the measure inequality which gives the required estimate.

We now turn to the key ideas in the proof of Theorem~\ref{thm:gt-linear}. Let \(L_1,L_2\) be the two omitted \(J\)-lines, with corresponding logarithmic potentials \(u_1,u_2\). Applying the technique from the proof of Theorem~\ref{thm:main}, together with a harmonic approximation for the relative potential $v=u_1-u_2$, we obtain the preliminary bound
\[
T_f(r)=O(r\log r),\qquad\text{and }\qquad a_f(r)=O(r\log r).
\]
To sharpen this growth rate,  we must therefore remove the logarithm from the relative potential. This is firstly achieved by a weighted Hodge--Calder\'on--Zygmund decomposition (see Lemma \ref{lem:weighted-splitting}), which yields
\[
M_r(|u_1-u_2|)=O(r).
\]
Since $P=L_1\cap L_2$ is omitted, the above divergence term $|F|$ can be bounded by $\|\mathrm{d}f\|$, providing a key tool to eliminate the log-error term in pair comparison.
Applying similar strategies in Lemma \ref{lem:perturbed}, the local estimate becomes genuinely linear, and we obtain
\[
M_r(u)=M_r\bigl(\max\{u_1,u_2\}\bigr)+O(r).
\]
Combining above together, it proves Theorem~\ref{thm:gt-linear}.

\medskip
\noindent\emph{Organization of the paper.}
Section~\ref{sec:value} recalls the value distribution theory in almost complex setting and constructs the $J$-line potentials.  Section~\ref{sec:pair} establishes the perturbed Da Costa--Duval potential lemma and the $O(r\log r)$ pair comparison.  Section~\ref{sec:main-proof} proves Theorem~\ref{thm:main}.  Finally, Section~\ref{sec:gt-proof} is devoted to Theorem~\ref{thm:gt-linear}: it first derives the preliminary logarithmic growth bound, then proves the localized linear pair comparison, performs the weighted splitting of the relative potential, and concludes the linear energy estimate.

\medskip
\noindent\emph{Notions and conventions.} Throughout this paper, we use the notion $\D(a,r)$ to denote the disc centered at $a$ with radius $r$. When $a=0$, we write $\D_r$ instead of $\D(0,r)$. The unit disc is denoted by $\D$. If $D=\D(a,r)$, then for each positive $\alpha>0$, by $\alpha D$ we mean the disc $\D(a,\alpha r)$.

\section{Value distribution theory and line potentials for \texorpdfstring{$J$}{J}-Brody curves}
\label{sec:value}
An almost complex structure $J$ on a real smooth manifold $X$ of real dimension $2n$  is a section of the bundle $\End(TM)$ satisfying $J^2=-\Id_{TM}$. The standard complex structure on $\mathbb{C}^n$ will be denoted by $i$. An almost complex structure on $\mathbb{CP}^2$ is said to be tamed by the Fubini-Study metric $\omega_{\FS}$ if $\omega_{\FS}(v,J(x)v)>0$ for any $v\in T_x(\mathbb{CP}^2)\setminus\{0\}$. A map $f\colon\mathbb{C}\rightarrow\mathbb{CP}^2$ of class $C^1$ is said to be an entire $J$-holomorphic curve if it satisfies $\mathrm{d} f\circ i=J\circ \mathrm{d}f$. Such an entire $J$-holomorphic curve $f$ is called a $J$-Brody curve if the norm of its derivative $\|\mathrm{d}f\|$ is uniformly bounded. By a $J$-lines in $\mathbb{CP}^2$, we mean a \(2\)-real-dimensional \(J\)-holomorphic submanifold representing the generator of \(H^2(\mathbb{CP}^2,\mathbb{Z})\).

In this section, we briefly recall the basis setup of value distribution theory for holomorphic $J$ curves and families of $J$-lines in the projective plane $(\mathbb{CP}^2,J)$. In particular, for a \(J\)-Brody curve \(f\colon \mathbb{C}\to(\mathbb{CP}^2,J)\), we restate its counting function with respect to a \(J\)-line \(L\)  according to the line potential, and we establish the First Main Theorem (FMT) with a sharper error term, corresponding to \cite{Duval-GAFA}. 

Throughout the paper we write
\[
\mathrm{d}_J^c\phi=-\frac1{2\pi}\mathrm{d}\phi\circ J,
\]
for any function \(\phi\) on \(\mathbb{CP}^2\). On the complex plane \(\mathbb{C}\), whose almost complex structure is always the standard one \(i\), we follow the usual notation that
\[
\mathrm{d}^c(\bullet)=\frac{i}{2\pi}(\bar\partial-\partial)(\bullet)
\qquad\text{and}\qquad
\mathrm{d}\mathrm{d}^c(\bullet)=\frac{i}{\pi}\partial\bar\partial(\bullet)
=\frac1{2\pi}\Delta(\bullet)\,\mathrm{d}x\wedge\mathrm{d}y,
\]
where \(\Delta\) is the Laplacian operator and \(z=x+iy\) is the affine coordinate on \(\mathbb{C}\).

For any real $2$-form $\alpha\in \mathcal{A}^{2}(\mathbb{C})$ on $\mathbb{C}$, the order function of \(\alpha\) is defined as
\[
T_\alpha(r)=\int_1^r\left(\int_{\D_t}\alpha\right)\frac{\mathrm{d}t}{t}.
\] 
For a \(J\)-Brody curve \(f\colon \mathbb{C}\to(\mathbb{CP}^2,J)\), recalling that its  Cartan characteristic function is given by
\[
T_f(r)=T_{f^*\omega_{\FS}}(r),
\]
where \(\omega_{\FS}\) is the Fubini-Study form on \(\mathbb{CP}^2\). Since $\mathbb{C}$ is contractible, we can write \(\mathrm{d}\mathrm{d}^cu=f^*\omega_{\FS}\) for certain function $u$ on $\mathbb{C}$. The Brody condition implies that $\Delta u$ is bounded. For a circle
$\Gamma=\partial \mathbb{D}(c,R)$ and an integrable function $\phi$ on $\Gamma$, we write its circular mean as
\[
M_{\Gamma}(\phi)
:=\frac{1}{2\pi R}\int_{\Gamma}\phi\,|\mathrm{d}z|.
\]
When $\Gamma=\Gamma_r:=\partial D(0,r)$ we abbreviate
$M_{\Gamma_r}(\phi)$ to $M_r(\phi)$.  Jensen's formula gives
\begin{equation}\label{eq:jensen-u}
	T_f(r)=M_r(u)+O(1).
\end{equation}

Let $\mathcal O(1)$ be the dual of the tautological line bundle $\mathcal O(-1)=\{(l,v)\in \mathbb{CP}^2\times\mathbb{C}\, |\, v\in l\}$, whose first Chern class is \(c_1(\mathcal{O}(1))=[\omega_{\FS}]\in H^2(\mathbb{CP}^2,\mathbb{Z})\). Note that it is a Hermitian complex line bundle, but not \(J\)-holomorphic. Fix a smooth Hermitian metric $h$ on $\mathcal O(1)$. For any $J$-line $L\subset\mathbb{CP}^2$, the Poincar\'e dual of \(L\) is \(\text{PD}([L])=c_1(\mathcal{O}(1))\). Then there exists a smooth section $\sigma$ of $\mathcal O(1)$ vanishing transversely on $L$; explicitly, $\|\sigma\|_h$ vanishes to first order along $L$. Without loss of generality, we assume $\|\sigma\|_h\le1$ and put
\[
\lambda=\log\|\sigma\|_h.
\]
Consequently, near \(L\) one has
\[\lambda(x)=\log\dist(x,L)+O(1).\]
Since \(\mathbb{CP}^2\) is compact, such \(O(1)\) term is uniform. Furthermore, the almost complex Poincar\'e-Lelong formula \cite[Section 4]{Duval-GAFA} gives that
\begin{equation}\label{eq:PL-line-target}
\mathrm{d}\mathrm{d}_J^c\lambda=[L]-\omega_{\FS}+\mathrm{d}\beta,
\end{equation}
where $\beta$ is a bounded smooth $1$-form away from $L$. In particular, $\beta$ is uniformly bounded in a tubular neighborhood of \(L\), hence uniformly bounded on $\mathbb{CP}^2\setminus L$ \cite[p. 755]{Duval-GAFA}.

Let \(\{L_i\}_{i=1,2,3}\) be three $J$-lines in general position and let \(\lambda_i,\beta_i\) be as above. Set \(u_i=u+\lambda_i\circ f\) for $i=1,2,3$. Then it follows from general position and compactness that
\begin{equation}\label{eq:triple-max}
u=\max\limits_{i=1,2,3}\{u_i\}+O(1).
\end{equation}
Pulling back \eqref{eq:PL-line-target} via the \(J\)-Brody curve \(f\), one has
\begin{equation}\label{eq:PL-line-source}
\mathrm{d}\mathrm{d}^cu_i=f^*[L_i]+\mathrm{d}\beta_{i,f},
\end{equation}
where \(\beta_{i,f}:=f^*\beta_i\). Since \(\beta_i\) is uniformly bounded on \(\mathbb{CP}^2\setminus L_i\) and \(f\) is Brody, the pullback \(\beta_{i,f}\) is uniformly bounded on \(\mathbb{C}\setminus f^{-1}(L_i)\). Consequently, for every $\mathbb{D}_t$ with \(t>0\), excising small disks around \(f^{-1}(L_i)\cap\D_t\) and using this uniform bound, the boundary terms tend to \(0\) as the radii shrink. And Stokes' theorem yields
\begin{equation}\label{eq:beta-stokes}
\left|T_{\mathrm{d}\beta_{i,f}}(r)\right|
\le \int_1^r
\left|\int_{\partial\D_t}\beta_{i,f}\right|\frac{\mathrm{d}t}{t}
=O(r).
\end{equation}
Define
\[
n_f(r,L_i)\coloneq f^*[L_i](\D_r)
\qquad\text{and}\qquad
N_f(r,L_i)\coloneq T_{f^*[L_i]}=\int_1^r n_f(t,L_i)\frac{\mathrm{d}t}{t}.
\]
Combining \eqref{eq:PL-line-source} and \eqref{eq:beta-stokes} implies that
\begin{align}\label{eq:first-main-line}
N_f(r,L_i)&=M_r(u_i)+O(r)\nonumber\\
&=T_f(r)+M_r(\lambda_i\circ f)+O(r),
\end{align}
which establishes FMT for \(J\)-Brody curve.
\begin{rmk}
In \cite[Section 6]{Duval-GAFA}, Duval proves FMT for \(J\)-holomorphic curve in the general elliptic setting, which is stated as \(N_f(r,L)=T_f(r)+M_r(\lambda\circ f)+o(T_f(r))\). His error term \(o(T)\) comes from the bounded \(1\)-form \(\beta\) and from the comparison of the characteristic functions \(T\) and \(T'\) via the Pl\"ucker estimates.
\end{rmk}
\begin{rmk}
	The relations \eqref{eq:PL-line-target} and \eqref{eq:PL-line-source} exhibit the main difference between the almost complex setting and the standard complex setting. When \(J=J_0\), the \(1\)-forms \(\beta\) and \(\beta_{i,f}\) vanish; consequently, \(u_i\) is subharmonic. For a general almost complex structure \(J\), the main technical difficulty arises from \(\beta\) and \(\beta_{i,f}\), as they produce a divergence perturbation.
\end{rmk}

\section{Pair potentials and the perturbed potential lemma}\label{sec:pair}
In this section, we prove some technical lemmas, which will be helpful in the sequel. For two distinct $J$-lines $L_i,L_j\subset\mathbb{CP}^2$ intersecting at the point \(P_{ij}=L_i\cap L_j\), we define
\[
\widehat\lambda_{ij}
=\frac12\log\bigl(\|\sigma_i\|_h^2+\|\sigma_j\|_h^2\bigr),
\]
admitting a logarithmic singularity at \(P_{ij}\). Recall that in \cite[pp.~2361]{Duval-Fourier04}, Duval constructs the central-projection
\[
\pi_{ij}:\mathbb{CP}^2\setminus\{P_{ij}\}\longrightarrow\mathbb P^1,
\]
parametrizing the pencil of $J$-lines through \(P_{ij}\) via the tangents. He also establishes the \(\mathrm{d}\mathrm{d}^c_J\)-formula that 
\begin{equation}\label{eq:central-PL}
\pi_{ij}^*\omega_{\mathbb P^1}
=\omega_{\FS}+\mathrm{d}\mathrm{d}_J^c\psi_{ij}+\mathrm{d}a_{ij}
\quad\text{on }\mathbb{CP}^2\setminus\{P_{ij}\},
\end{equation}
where \(\psi_{ij}(x)\) is a negative function identifying with \(\widehat\lambda_{ij}(x)\) around \(P_{ij}\), and $a_{ij}$ is a uniformly bounded $1$-form on \(\mathbb{CP}^2\setminus\{P_{ij}\}\). Moreover, one has \(\psi_{ij}=\widehat\lambda_{ij}+O(1)\) on \(\mathbb{CP}^2\setminus\{P_{ij}\}\).
Set
\begin{equation}\label{eq:defi-uij}
	u_{ij}=u+\psi_{ij}\circ f.
\end{equation}
Since \(\widehat\lambda_{ij}=\max\{\lambda_i,\lambda_j\}+O(1)\), we have
\begin{equation}\label{eq:uij-max}
u_{ij}=\max\{u_i,u_j\}+O(1).
\end{equation}

\begin{lem}\label{lem:pair-potential}
	There exist a positive \((1,1)\)-current $\mu_{ij}$ with locally finite mass and a bounded $1$-form $\alpha_{f,ij}$ on $\C$ such that
	\begin{equation}\label{eq:uij-decomp}
	\mathrm{d}\mathrm{d}^cu_{ij}=\mu_{ij}+\mathrm{d}\alpha_{f,ij}.
	\end{equation}
\end{lem}
\begin{proof}
	On $\C\setminus f^{-1}(P_{ij})$, combining \eqref{eq:central-PL} and \eqref{eq:defi-uij}, one has
	\[
		\mathrm{d}\mathrm{d}^cu_{ij}=(\pi_{ij}\circ f)^*\omega_{\mathbb P^1}-f^*\mathrm{d}a_{ij}
		\]
	Write \(\alpha_{f,ij}:=-f^*a_{ij}\) defined on $\C\setminus f^{-1}(P_{ij})$, hence the Brody condition and boundedness of $a_{ij}$ give that
	$\|\alpha_{f,ij}\|_{L^\infty}\le C$.  Since $\pi_{ij}\circ f$ is quasiconformal \cite[pp. 2362]{Duval-Fourier04}, the \(2\)-form $(\pi_{ij}\circ f)^*\omega_{\mathbb{CP}^1}$ is a nonnegative area form (with the chosen orientation). Hence we define \(\mu_{ij}\coloneq (\pi_{ij}\circ f)^*\omega_{\mathbb P^1}\) and obtain \eqref{eq:uij-decomp} on the punctured domain $\C\setminus f^{-1}(P_{ij})$.
	
	To control the singularity at $f^{-1}(P_{ij})$, we consider the almost complex blow-up $\sigma_{ij}\colon X_{ij}\rightarrow\mathbb{CP}^2$ at the point $P_{ij}$, and denote the exceptional divisor $E_{ij}=\sigma_{ij}^{-1}(P_{ij})$. Let $\tilde{f}\colon\mathbb{C}\rightarrow X_{ij}$ be the lifting of $f$, and let $\tilde{J}$ be the lifted almost complex structure in \(X_{ij}\), admitting a Lipschitz extension along \(E_{ij}\). Since $\pi_{ij}$ extends to a global well-defined smooth fibration $\tilde{\pi}_{ij}\coloneq \pi_{ij}\circ \sigma_{ij}\colon X_{ij}\rightarrow\mathbb{CP}^1$, the composition $\tilde{\pi}_{ij}\circ \tilde{f}\colon\mathbb{C}\rightarrow\mathbb{CP}^1$ is a well-defined smooth map \cite[pp. 2361]{Duval-Fourier04}.
	\begin{claim}
	On the blowup \(X_{ij}\), one has
	\begin{equation}
		\label{saleur-pullback}
		\tilde{\pi}_{ij}^*\omega_{\mathbb{CP}^1} = \sigma_{ij}^*\omega_{\FS} + \mathrm{d}\mathrm{d}^c_{\tilde{J}}(\sigma_{ij}^*\psi_{ij}) - [E_{ij}] + \mathrm{d}(\sigma_{ij}^*{a}_{ij}).
		\end{equation}
	\end{claim}
	\begin{proof}[Proof of Claim]
	Pulling back \eqref{eq:central-PL} to $X_{ij}\setminus E_{ij}$, where $\sigma_{ij}$ is a diffeomorphism onto $\mathbb{CP}^2\setminus\{P_{ij}\}$, it gives the identity of smooth forms
	\begin{equation}
	\label{eq:pullback-complement}
	\tilde\pi_{ij}^*\omega_{\mathbb{CP}^1}
	=
	\sigma_{ij}^*\omega_{\FS}
	+
	\mathrm{d}\mathrm{d}_{\tilde J}^c(\sigma_{ij}^*\psi_{ij})
	+
	\mathrm{d}(\sigma_{ij}^*a_{ij})
	\quad\text{on }X_{ij}\setminus E_{ij}.
	\end{equation}
	We then extend it across the exceptional divisor $E_{ij}$, with last two terms potentially singular.
	
	We firstly consider \(\mathrm{d}(\sigma_{ij}^*a_{ij})\). Via a diffeomorphism from \((\mathbb{CP}^2,P_{ij})\) to \((\mathbb{C}^2,0)\), it suffices to consider the complex blowup of \(\mathbb{C}^2\) around \(0\). According to \cite[p. 2361]{Duval-Fourier04}, around \(P_{ij}\) equipped with local coordinate \((x,w)\), we have the local form \(a_{ij}(x,w)=\mathrm{d}\psi_{ij}\circ O(|(x,w)|)+b(x,w)\) with $b(x,w)$ bounded. Hence in the blow-up chart \((u,v)\) such that $x=uv$, $w=v$, we have
	\[
	\sigma_{ij}^*a_{ij}	=\Bigl(\frac{\mathrm{d}v}{v}+\omega(u)\Bigr)\circ O(|v|)+b(uv,v),
	\]
	where \(\omega(u)\) is a smooth \(1\)-form. The first term is bounded as $v\to 0$, and it vanishes on $E_{ij}$. The second term $b(uv,v)$ is bounded anywhere. Hence $\sigma_{ij}^*a_{ij}$ extends to a bounded $1$-form on the whole $X_{ij}$, therefore \(\mathrm{d}(\sigma_{ij}^*a_{ij})\) is well-defined on $X_{ij}$.
	
	The function $\sigma_{ij}^*\psi_{ij}$ has a logarithmic singularity along $E_{ij}$. A local computation in the blow-up charts and Poincar\'e–Lelong formula show that
	\begin{equation}
	\label{eq:PL-blowup}
	\mathrm{d}\mathrm{d}_{\tilde J}^c(\sigma_{ij}^*\psi_{ij})=[E_{ij}]+\omega_{\text{smooth}},
	\end{equation}
	where $\omega_{\text{smooth}}$ is a smooth $(1,1)$-form. Note that the singular part $[E_{ij}]$ is independent of the almost complex structure $\tilde J$. Since the other terms in \eqref{eq:pullback-complement} are all smooth across $E_{ij}$, comparing \eqref{eq:pullback-complement} and \eqref{eq:PL-blowup} on the dense open set $X_{ij}\setminus E_{ij}$, it yields \eqref{saleur-pullback} on the whole \(X_{ij}\).
	\end{proof}
	  
	Pulling back \eqref{saleur-pullback} again via $\tilde{f}$ and combining with \eqref{eq:defi-uij}, we obtain
	\begin{align*}
	(\tilde{\pi}_{ij} \circ \tilde{f})^*\omega_{\mathbb{CP}^1}
	&= 
	\tilde{f}^*(\sigma_{ij}^*\omega_{\FS})	+\tilde{f}^*\mathrm{d}\mathrm{d}^c(\sigma_{ij}^*\psi_{ij})- \tilde{f}^*[E_{ij}]+\tilde{f}^*(\mathrm{d}\sigma_{ij}^*{a}_{ij})\nonumber\\
	&=(\sigma_{ij}\circ\tilde{f})^*\omega_{\FS}+\mathrm{d}\mathrm{d}^c\big((\sigma_{ij}\circ\tilde{f})^*\psi_{ij}\big)	-\tilde{f}^*[E_{ij}]+\mathrm{d}(\tilde{f}^*(\sigma_{ij}^*{a}_{ij}))\nonumber\\
	&=\mathrm{d}\mathrm{d}^cu+\mathrm{d}\mathrm{d}^c(\psi_{ij}\circ f)- \tilde{f}^*[E_{ij}]+\mathrm{d}(\tilde{f}^*(\sigma_{ij}^*{a}_{ij}))\nonumber\\
	&=\mathrm{d}\mathrm{d}^cu_{ij}-\tilde{f}^*[E_{ij}]+\mathrm{d}(\tilde{f}^*(\sigma_{ij}^*{a}_{ij})).
	\end{align*}
	The current $\tilde{f}^*[E_{ij}]=\sum\limits_{z\in f^{-1}(P_{ij})}	\mu_z\delta_z$ is a sum of Dirac masses, where $\mu_z\geq1$ is the local multiplicity of \(f\) at $z$. Putting $\alpha_{f,ij}=\tilde{f}^*(\sigma_{ij}^*{a}_{ij})$, we obtain at the end the desired equality
	\[
	\mathrm{d}\mathrm{d}^cu_{ij}
	=(\tilde{\pi}_{ij} \circ\tilde{f})^*\omega_{\mathbb{CP}^1}+\sum_{z\in f^{-1}(P_{ij})}	\mu_z\delta_z+\mathrm{d}\tilde{f}^*(\sigma_{ij}^*{a}_{ij})\eqcolon \mu_{ij}+\mathrm{d}\alpha_{f,ij}.
	\]
	Indeed, the boundedness of \(\alpha_{f,ij}\) can be justified as follows: we consider locally in the blow-up chart \((u,v)\) with \(x=uv\) and \(w=v\) as above. Pulling back by \(\tilde{f}=(u(z),v(z))\) yields
	\[
	\tilde{f}^*(\sigma_{ij}^*a_{ij})=O(|v|)\left(\frac{v'(z)\mathrm{d}z}{v(z)}+u'(z)\mathrm{d}z\right)+O(1).
	\]
Since \(f\) is Brody, we have \(|v'|=|w'|\le C\), and \(|u'v| \le C\) since \(|x'|=|u'v+uv'|\le C\).
\end{proof}

We next prove the local estimate responsible for the $O(r\log r)$ error term, which is an analog of \cite[Lemma 2]{Dacosta-Duval2013}. Recall that the circular mean \(M_{\Gamma}(\phi)=\frac{1}{2\pi R}\int_{\Gamma}\phi\,|\mathrm{d}z|\), where $\Gamma=\partial \mathbb{D}(c,R)$ and $\phi$ is an integrable function on $\Gamma$. In the following context, \(C\) represents different constants.
\begin{lem}\label{lem:perturbed}
	Let $D=\mathbb{D}(a,\rho)$ and $2D=\mathbb{D}(a,2\rho)$. Let $z_0\in\partial(2D)$. Suppose that $w$ is a non-negative function on $2D$, $w(z_0)=0$ and $w$ is of $\mathcal{C}^1$-class near $z_0$. Assume that, in the sense of distributions, one has
	\begin{equation}\label{eq:perturbed-PDE}
		-\Delta w=\mu+\operatorname{Div}F,
	\end{equation}
	where $\mu$ is a positive Radon measure and $\|F\|_{L^\infty(2D)}\le A$. Assume also that for some $0<r_0\le\rho$, it holds
	\begin{equation}\label{eq:F-local-lip}
		|F(z)-F(z_0)|\le L|z-z_0|
		\quad\text{on }2D\cap D(z_0,r_0).
	\end{equation}
	Then for every arc $\gamma=\Gamma\cap D$ whose containing circle has radius $R\ge\rho$, there exists a universal constant \(C\) satisfying
	\begin{equation}\label{eq:perturbed-general}
		M_\gamma(w)\coloneq\frac{1}{2\pi R}\int_{\gamma}w\,|\mathrm{d}z|
		\le C\frac{\rho^2}{R}\left(
		|\nabla w(z_0)|+A+Lr_0+A\log\frac{4\rho}{r_0}
		\right).
	\end{equation}
	
Furthermore, if \eqref{eq:F-local-lip} is available with uniform constants on a fixed radius $r_*>0$ (independent of $\rho$), then with $r_0=\min\{r_*,\rho\}$, we have
	\begin{equation}\label{eq:perturbed-log}
		M_\gamma(w)\le C\frac{\rho^2}{R}\bigl(1+\log(2+\rho)\bigr).
	\end{equation}
	The same conclusion holds if $w$ has isolated positive logarithmic poles in $2D$.
\end{lem}

\begin{proof}
	By Green's representation Theorem, one has 
	\[
	w(x)=\int_{\partial(2D)} w(t)\partial_\nu G(x,t)\, \mathrm{d}s(t)-\int_{2D}G(x,t) \Delta_t w(t) \, \mathrm{d}t,
	\]
	where \(\nu\) denotes the inward unit normal to \(\partial(2D)\), and $G(x,t)$ is the Green function of $2D$, normalized by
	\[
	-\Delta_tG(x,t)=\delta_x,
	\qquad	G(x,t)=0\quad\forall\, t\in\partial(2D).
	\]
	Let $H$ be the harmonic function on $2D$ with the boundary values of $w$. Substituting \eqref{eq:perturbed-PDE}, it yields that
	\begin{align}\label{eq:green-rep}
		w(x)&=H(x)+\int_{2D}G(x,t)\,\mathrm{d}\mu(t)+\int_{2D} G(x,t)\operatorname{Div} F(t)\,\mathrm{d}t\nonumber\\
		&=H(x)+\int_{2D}G(x,t)\,\mathrm{d}\mu(t)-\int_{2D} \nabla_tG(x,t)\cdot F(t)\,\mathrm{d}t.
	\end{align}
	For every \(x \in 2D\), using the divergence theorem on the constant vector field \(F(z_0)\), we have
		\[
		\int_{2D}\nabla_t G(x,t)\cdot F(z_0)\,\mathrm{d}t=\int_{\partial(2D)}G(x,t)(F(z_0)\cdot\nu)\,\mathrm{d}s(t) = 0.
		\]
	Therefore, we can replace \(F(t)\) by \(\tilde{F}(t)=F(t)-F(z_0)\), i.e.,
	\begin{equation}
		\label{Green representation theorem, replace F by tilde F}
		w(x) = H(x) + \int_{2D} G(x,t) \, \mathrm{d}\mu(t) - \int_{2D} \nabla_t G(x,t) \cdot \tilde{F}(t) \, \mathrm{d}t.
	\end{equation}
	
	Denote by \(P(t)=\partial_{\nu(z_0)}G(z_0,t)=\frac{1}{4\pi\rho}\,\frac{(2\rho)^2-|t-a|^2}{|t-z_0|^2}\), the Poisson kernel on \(2D\) valued at $z_0$. Since $\tilde{F}$ vanishes at \(z_0\) and satisfies the Lipschitz estimate \eqref{eq:F-local-lip}, the product $\nabla P(t)\cdot \tilde{F}(t)$ is locally integrable (see below). Hence, we can differentiate under the integral sign in the inward normal direction \(\nu(z_0)\), deriving that
	\begin{equation*}
		\partial_\nu w(z_0)
		=\partial_\nu H(z_0)+\int_{2D}P(t)\,\mathrm{d}\mu(t)
		-\int_{2D}\nabla P(t)\cdot\widetilde F(t)\,\mathrm{d}t.
	\end{equation*}
	Since $H\ge0$, $H(z_0)=0$, and $\mu\ge0$, the first two terms on the right are non-negative. Also $w\ge0$ and $w(z_0)=0$, we then have $0\le\partial_\nu w(z_0)\le|\nabla w(z_0)|$. Hence
	\begin{equation}\label{eq:positive-normal-control}
		\partial_\nu H(z_0)+\int_{2D}P\,\mathrm{d}\mu
		\le |\nabla w(z_0)|+\left|\int_{2D}\nabla P(t)\cdot\widetilde F(t)\,\mathrm{d}t\right|.
	\end{equation}

	Firstly, one has \(|\nabla P(t)|\le \frac{C}{|t-z_0|^2}.\)
	Combined with \eqref{eq:F-local-lip} on the region $2D\cap \mathbb{D}(z_0,r_0)$, it gives that
	\[
	\int_{2D\cap \mathbb{D}(z_0,r_0)}|\nabla P|\,|\widetilde F|\, \mathrm{d}t
	\le CL\int_0^{r_0}\frac{r}{r^2}\,r\,\mathrm{d}r
	\le CLr_0.
	\]
	On the complement, $|\widetilde F|\le2A$, hence we have
	\[
	\int_{2D\setminus \mathbb{D}(z_0,r_0)}|\nabla P|\,|\widetilde F|\, \mathrm{d}t
	\le CA\int_{r_0}^{4\rho}\frac{\mathrm{d}r}{r}
	=CA\log\frac{4\rho}{r_0}.
	\]
	Consequently
	\begin{equation}\label{eq:E-bound}
		\left|\int_{2D}\nabla P(t)\cdot\widetilde F(t)\,\mathrm{d}t\right|\le C\left(Lr_0+A\log\frac{4\rho}{r_0}\right).
	\end{equation}
	
	For every arc $\gamma=\Gamma\cap D$ whose containing circle has radius $R\ge\rho$, we define
	\[
	K_\gamma(t)=\frac1{2\pi R}\int_\gamma G(x,t)\,|\mathrm{d}x|.
	\]
	If \(t\in\tfrac32 D\) and \(x\in D\), the disk Green kernel satisfies \(G(x,t)\le C\bigl(1+\log^+(\frac{\rho}{|x-t|})\bigr)\). Moreover, \(\operatorname{Length}(\gamma)\le C\rho\), and for \(0<s\le\rho\) the arc length satisfies \(\operatorname{Length}(\Gamma\cap \mathbb{D}(t,s))\le Cs\). Using the formula \(\log^+(\frac{\rho}{|x-t|})= \int_{|x-t|}^{\rho}\frac{\mathrm{d}s}{s}\), we integrate over \(\gamma\) to get
	\[
	\int_\gamma\log^+\frac{\rho}{|x-t|}\,|\mathrm{d}x|
	= \int_0^\rho \frac{\operatorname{Length}(\gamma\cap \mathbb{D}(t,s))}{s}\,\mathrm{d}s
	\le \int_0^\rho \frac{Cs}{s}\,\mathrm{d}s
	\le C\rho,
	\]
	which yields \(K_\gamma(t)\le C\frac{\rho}{R}\) for \(t\in\tfrac32 D\). On $\tfrac32 D$, the Poisson kernel satisfies \(P(t)=\frac{1}{4\pi\rho} \frac{(2\rho)^2-|t-a|^2}{|t-z_0|^2}\ge \frac{c}{\rho}\), and therefore \(K_\gamma(t)\le C\frac{\rho^2}{R}P(t)\). To extend this inequality to all of \(2D\), we consider the function 
	\[v(t) = C\frac{\rho^2}{R}P(t) - K_\gamma(t).\] 
	Both \(K_\gamma\) and \(P\) are harmonic on \(2D\setminus\tfrac32 D\), and they vanish on \(\partial(2D)\setminus\{z_0\}\). On the inner boundary \(\partial(\tfrac32 D)\), we have \(v(t)\ge 0\) by the estimate above. On the outer boundary \(\partial(2D)\setminus\{z_0\}\), since \(K_\gamma=0\) and \(P\ge0\), we also have \(v(t)\ge0\). Although \(P(t)\) has a positive boundary singularity at \(z_0\), the function \(v(t)\) cannot attain a negative minimum there. Therefore, by the maximum principle, \(v(t)\ge0\) on \(2D\setminus\tfrac32 D\), which extends the bound
	\[0\le K_\gamma(t)\le C_1(\frac{\rho^2}{R})P(t),\quad \forall\, t\in 2D.\]
	
	For \(x\in D\), the Green kernel satisfies \(|\nabla_tG(x,t)|\le \frac{C}{|x-t|}\),
	which implies
	\begin{align}\label{eq:green-kernel-estimate}
	\int_{2D} |\nabla_t G(x,t)| \, \mathrm{d}t \le C \int_{2D} \frac{\mathrm{d}t}{|x-t|}
	\leq C  \int_{\mathbb{D}(x, 3\rho)} \frac{\mathrm{d}t}{|x-t|}=6C\pi\rho.
	\end{align}
	Therefore we have
	\[
	\begin{aligned}
		\int_{2D} |\nabla K_\gamma(t)| \, \mathrm{d}t \le \frac{1}{2\pi R} \int_\gamma \left( \int_{2D} |\nabla_t G(x,t)| \, \mathrm{d}t \right) |\mathrm{d}x| \le \frac{C \rho}{R} \cdot \operatorname{Length}(\gamma)\le C\frac{\rho^2}{R}.
	\end{aligned}
	\]
	
	Averaging \eqref{eq:green-rep} over $\gamma$, it gives
	\begin{equation*}
		M_\gamma(w)=M_\gamma(H)+\int_{2D}K_\gamma\,\mathrm{d}\mu
		-\int_{2D}\nabla K_\gamma\cdot {F}\,\mathrm{d}t.
	\end{equation*}
	By Harnack inequality for the positive harmonic function $H$, we have \(\sup_DH\le C\rho\,\partial_\nu H(z_0)\). Thus
	\[
	M_\gamma(H)
	\le \frac{\operatorname{Length}(\gamma)}{2\pi R}\sup_D H
	\le C\frac{\rho^2}{R}\,\partial_\nu H(z_0),
	\]
and
	\[
	\int_{2D}K_\gamma\,\mathrm{d}\mu\le C\frac{\rho^2}{R}\int_{2D}P\,\mathrm{d}\mu,
	\qquad
	\left|\int_{2D}\nabla K_\gamma\cdot {F}\,\mathrm{d}t\right|\le CA\frac{\rho^2}{R}.
	\]
	Combining above with \eqref{eq:positive-normal-control} and \eqref{eq:E-bound}, we proves \eqref{eq:perturbed-general}. The fixed-radius form \eqref{eq:perturbed-log} directly follows by taking $r_0=\min\{r_*,\rho\}$.

	Finally, we consider the case where \(w\) has isolated positive logarithmic poles in \(2D\). Suppose near the pole \(p \in 2D\), \(w(z) = c\,\log(\frac{1}{|z-p|}) + O(1)\) with \(c>0\).  Let \(2D_\varepsilon=2D\setminus\overline{\mathbb{D}(p,\varepsilon)}\) be the punctured domain, and let \(G_\varepsilon(x,t)\) be the Green function on \(2D_\varepsilon\). Then Green representation theorem yields
	\begin{align*}
	w(x) =& \int_{\partial(2D)} w(t)\partial_\nu G_\varepsilon(x,t)\, \mathrm{d}s(t) 
	+ \int_{\partial \mathbb{D}(p,\varepsilon)} w(t)\partial_\nu G_\varepsilon(x,t)\, \mathrm{d}s(t) \\
	&- \int_{2D_\varepsilon}G_\varepsilon(x,t) \Delta_t w(t) \, \mathrm{d}t.
	\end{align*}
	As \(\varepsilon \to 0\), the second term converges to \(0\). In the sense of distributions, the log- singularity at \(p\) contributes a positive Dirac mass. Therefore, the third term converges to
	\[\int_{2D} G(x,t) \bigl(2\pi c\,\delta_p+\mu+\text{Div}F \bigr) \, \mathrm{d}t.\]
	The Dirac mass is locally integrable and can be absorbed by the positive measure \(\mu\). Therefore \eqref{eq:green-rep} still holds and the subsequent proof remains valid.
\end{proof}

Now we prove the main result of this section.
\begin{pro}\label{prop:two-potential}
	For any two distinct $J$-lines $L_i,L_j$ among a fixed finite family in general position,
	\begin{equation*}
		M_R(u)=M_R\bigl(\max\{u_i,u_j\}\bigr)+O\bigl(R\,\log(2+R)\bigr).
	\end{equation*}
\end{pro}
\begin{proof}
	Let $L_k$ be a $J$-line in general position with $L_i,L_j$. 
	It follows from \eqref{eq:triple-max} and \eqref{eq:uij-max} that,
	\[
	u=\max\{u_{ij},u_k\}+O(1)=u_{ij}+(u_k-u_{ij})_++O(1).
	\]
	Rescale the defining sections \(\sigma_i,\sigma_j,\sigma_k\) such that the set \(\Omega\coloneq \{z\in \mathbb{C}\,|\,u_k(z)>u_{ij}(z)\}\subsetneq\mathbb{C}\) is non-empty. And define \(w\coloneq u_k-u_{ij}\) as a positive function on \(\Omega\), whether or not \(f^{-1}(P_{ij})\in \Omega\) (i.e. whether or not \(w\) has a logarithmic pole).
	
	The set \(f^{-1}(L_k)\) does not meet $\Omega$ because \(u_k=-\infty\) there, whereas \(u_{ij}\) is finite by general position. Hence on $\Omega$, it follows from \eqref{eq:PL-line-source} and Lemma~\ref{lem:pair-potential} that
	\begin{equation}\label{eq:w-current}
	-\mathrm{d}\mathrm{d}^cw=-\mathrm{d}\mathrm{d}^c(u_k-u_{ij})=\mu_{ij}+\mathrm{d}\alpha_{f,ij}-\mathrm{d}f^*\beta_k,
	\end{equation}
	where \(f^*\beta_k,\alpha_{f,ij}\) are bounded \(1\)-forms. In real coordinates, after multiplying by the normalization constant, \eqref{eq:w-current} can be rewritten in the sense of distributions as 
	\begin{equation}\label{eq:w-scalar}
	-\Delta w=\widetilde\mu+\operatorname{Div}F,
	\end{equation}
	where \(\widetilde\mu\) is a positive measure and \(\|F\|_\infty\le A\), with the coefficients of \(F\) depending on \(\alpha_{f,ij}\) and \(f^*\beta_k\).
	
	Consider the set \(K=\{x\in\mathbb{CP}^2:\lambda_k(x)=\psi_{ij}(x)\}\) such that \(f(\partial\Omega)\subset K\). Observe that \(K\) is contained in a compact subset of \(\mathbb{CP}^2\setminus\bigl(L_k\cup\{P_{ij}\}\bigr)\). On a fixed smooth neighborhood of $K$, the function $\lambda_k-\psi_{ij}$ is smooth with uniformly bounded derivative. Combining with the Brody bound of \(f\) gives that, for \(z_0\in \partial\Omega\),
		\begin{equation}\label{eq:contact-gradient}
			|\nabla w(z_0)|
			=|\mathrm{d}(\lambda_k-\psi_{ij})_{f(z_0)}\circ \mathrm{d}f(z_0)|
			\le C.
		\end{equation}
	Choose $r_*>0$ sufficiently small such that every $f(\mathbb{D}(z_0,r_*))$ is contained in the fixed smooth neighborhood of \(K\). Since \(f\) is \(J\)-holomorphic, standard interior elliptic estimates \cite{McDuff91} yields a uniform $\mathcal{C}^2$-bound for $f$ on the disk \(\mathbb{D}(z_0,r_*)\). Hence we derive the Lipschitz estimate that, for any \(z\in \mathbb{D}(z_0,r_*)\),
	\begin{equation}\label{eq:F-contact-lip}
		|F(z)-F(z_0)|\le L|z-z_0|.
	\end{equation}
	
	Based on \eqref{eq:w-scalar}-\eqref{eq:F-contact-lip} we now apply Theorem \ref{lem:perturbed} as follows. Choose \(R\) appropriately large and translate \(\Gamma_R\) if necessary such that, for any $a\in\Gamma_R\cap\Omega$ and $\rho(a)=\frac{\dist(a,\partial\Omega)}{2}$, we have \(\rho(a)\leq R\). Define \(D_a \coloneq \mathbb{D}(a,\rho(a))\) and \(2D_a \coloneq \mathbb{D}(a,2\rho(a))\). Then \(2D_a\subset\Omega\) and there exists \(z_0\in\partial(2D_a)\cap\partial\Omega\) with \(w(z_0)=0\). 
	Set \(\gamma(a)=\Gamma_R\cap D_a\). Then Lemma~\ref{lem:perturbed} gives that
	\begin{equation}\label{eq:local-rlogr}
		M_{\gamma(a)}(w)
		\le C\frac{\rho(a)^2}{R}\bigl(1+\log(2+\rho(a))\bigr).
	\end{equation}
	Applying the Besicovitch covering theorem \cite{bes45} to the family of disks $\{D_a:a\in\Gamma_R\cap\Omega\}$ which covers the bounded set \(\Gamma_R\cap\Omega\), there exist finite subcollections of at most countable disjoint disks, whose union still covers \(\Gamma_R \cap \Omega\). Merging these subcollections, we obtain a countable subfamily \(\{D_\nu = \mathbb{D}(a_\nu, \rho_\nu)\}_\nu\) covering \(\Gamma_R \cap \Omega\) with uniformly finite overlap. Since \(\rho_\nu\leq R\) and each arc \(\gamma_\nu = \Gamma_R \cap D_\nu\) has length comparable from below to \(\rho_\nu\), we have
	\begin{equation*}\label{eq:radius-packing}
		\sum_\nu\rho_\nu\le CR.
	\end{equation*}
	Summing \eqref{eq:local-rlogr} by \(\nu\), we obtain
	\begin{align*}
		M_R((u_k-u_{ij})_+)=M_R(w_+)
		&\le \sum_\nu M_{\Gamma_R\cap D_\nu}(w)\\
		&\le \frac CR\sum_\nu\rho_\nu^2\bigl(1+\log(2+\rho_\nu)\bigr)\\
		&\le C\sum_\nu\rho_\nu\bigl(1+\log(2+\rho_\nu)\bigr)\\
		&\le CR\,\log(2+R),
	\end{align*}
Recall that \(u=u_{ij}+(u_k-u_{ij})_++O(1)\), then we finish the proof.
\end{proof}

\section{Proof of Theorem \ref{thm:main}}\label{sec:main-proof}
We follow \cite[Theorem 3]{Dacosta-Duval2013} to complete the proof of Theorem \ref{thm:main}. A function is called $\delta$-subharmonic if locally, it can be written as the difference of two subharmonic functions.

\begin{lem}[\cite{Grishin83,Eremenko-Sodin91}]\label{lem:grishin}
	Let $\varphi\ge0$ be a $\delta$-subharmonic function on a plane domain and let $E\subset\{\varphi=0\}$ be a Borel set. Then the restriction of the Riesz charge $\Delta\varphi$ to $E$ is a nonnegative measure.
\end{lem}

\begin{lem}[{\cite[Lemma~5]{Dacosta-Duval2013}.}]\label{lem:dcd-zero-set}
	Let $\varphi\geq 0$ be a lower semi-continuous $\delta$-subharmonic function on $2\D$ such that $\Delta\varphi$ is $L^{\infty}$ on $\{\varphi = 0\}$, then $\Delta\varphi = 0$ on the set $\{\varphi=0\}$.
\end{lem}
\begin{proof}[Proof of Theorem \ref{thm:main}]
	It suffices to prove the non-integrated form
	\begin{equation*}
	(q-1)a_f(r)\le \sum_{i=1}^q n_f(r,L_i)+o(r^2),
	\end{equation*}
	where \(a_f(r)=\int_{\D_r}f^*\omega_{\FS}\). Equivalently, we are supposed to show that	
	\begin{equation}\label{eq:limsup}
	\limsup_{r \to \infty} \frac{1}{r^2} \left[(q-1)a_f(r) - \sum_{i=1}^q n_f(r, L_i)\right] \le 0.
	\end{equation}
	Let $\{r_\nu\}_{\nu}\to\infty$ be a sequence of radii such that the LHS of \eqref{eq:limsup} converges. Recall that \(f^*\omega_{\FS}=\mathrm{dd}^cu\), and that in \eqref{eq:PL-line-source} each \(\beta_i\) is bounded for \(i=1,\cdots,q\). Since \(f\) is Brody, we obtain a converging sequence
	\[
	\dfrac{1}{r_{\nu}^2}\big[(q-1)\Delta u-\sum_{i=1}^q\Delta u_i\big](\D_{r_{\nu}}).
	\]
	Suppose that \(f(0) \notin \cup_{i=1}^qL_i\) (i.e. \(u_i(0)\) is finite for any \(i=1,\ldots,q\)). Let $h_\nu$ be the least harmonic majorant of $u$ on $\D_{2r_\nu}$ and put $\Lambda_\nu(z)=r_\nu z$. We then consider the sequences of functions on the disk  $2\D$
	\begin{equation}\label{eq:scaled-v}
	v_\nu(z)=\frac{u(r_\nu z)-h_\nu(r_\nu z)}{r_\nu^2},
	\qquad
	v_{i,\nu}(z)=\frac{u_i(r_\nu z)-h_\nu(r_\nu z)}{r_\nu^2}.
	\end{equation}
	
	Here each $v_\nu$ is a non-positive subharmonic function and $\Delta v_\nu$ is uniformly bounded, since $\Delta u$ is bounded. Hence by compactness argument~\cite{Hormander03}, either $\{v_\nu\}$  converges locally uniformly to $-\infty$, or after extracting a subsequence, it converges to a subharmonic function $v$ in the sense of $L^1_{\text{loc}}$. Since $f$ is Brody, $h_\nu(0)=\int_{\partial\D_{2r_\nu}}u(0)=O(r_\nu^2)$, hence $v_\nu(0)$ remains bounded as $\nu\to\infty$. Thus the first case is ruled out.
	
	By \eqref{eq:PL-line-source}, we have \(\mathrm{d}\mathrm{d}^cv_{i,\nu}=\frac1{r_\nu^2}\Lambda_\nu^*(f^*[L_i])+\mathrm{d}(\frac1{r_\nu^2}\Lambda_\nu^*\beta_{i,f})\), where \(\|\frac1{r_\nu^2}\Lambda_\nu^*\beta_{i,f}\|_\infty\le \frac{C}{r_\nu}\to 0\). Hence \(v_{i,\nu}\) is not subharmonic due to a small divergence perturbation \(\mathrm{d}(\frac1{r_\nu^2}\Lambda_\nu^*\beta_{i,f})\). To apply the subharmonic compactness theorem, we consider an arbitrary relatively compact subdomain $D'\Subset2\D$ and solve the Dirichlet problem as follows:
	\[
	-\Delta p_{i,\nu} = \operatorname{Div}(\frac1{r_\nu^2}\Lambda_\nu^*\beta_{i,f}) \quad \text{in } D', 
	\qquad p_{i,\nu} = 0 \quad \text{on } \partial D'.
	\]
	Therefore, the function \(w_{i,\nu} \coloneq v_{i,\nu}-p_{i,\nu}\) is subharmonic on \(D'\), with \(\|p_{i,\nu}\|_{L^\infty}\le C\|\frac1{r_\nu^2}\Lambda_\nu^*\beta_{i,f}\|_\infty\to 0\). Hence by \cite{Hormander03}, either \(w_{i,\nu} \to -\infty\), or it converges to a subharmonic function \(v_i\) in \(L^1_{\mathrm{loc}}\). Since \(f(0) \notin L_i\) and \(\|p_{i,\nu}\|_\infty \to 0\), we have \(w_{i,\nu}(0)=v_{i,\nu}(0)-p_{i,\nu}(0)= v_\nu(0)+O(r_\nu^{-2})-p_{i,\nu}(0)\), remaining bounded. Thus it rules out the first case and \(v_{i,\nu}\to v_i\) in \(L^1_{\mathrm{loc}}(D')\). By a standard diagonal extraction over an exhaustion of \(2\D\) by compact subdomains, it yields a subsequence such that \(v_{i,\nu} \to v_i\) in \(L^1_{\mathrm{loc}}(2\D)\).
	
	Since the defining sections were normalized by $\|\sigma_i\|_h\le1$, we have $u_i\le u$ and hence $v_{i,\nu}\le v_\nu$. Passing to the limit, it gives that $v_i\le v$. By \eqref{eq:triple-max}, we have \(v_\nu = \max\{v_{i,\nu}, v_{j,\nu}, v_{k,\nu}\} + O(r_\nu^{-2})\). Passing to the limit, we obtain $v =\max\{v_i, v_j, v_k\}$. We claim that
	\begin{equation}\label{eq:pair-limit}
	v=\max\limits_{i\ne j}\{v_i,v_j\}.
	\end{equation}
	To see this, fix an arbitrary point \(c\in 2\D\) and \(0<s_0<s_1\) such that the closed annulus
	\[
	A(c;s_0,s_1)=\{z:s_0\le|z-c|\le s_1\}\Subset 2\D.
	\]
	For \(t\in[s_0,s_1]\), applying Proposition~\ref{prop:two-potential} to the circle \(\partial \D(c,t)\), we receive
	\[
	0\le M_{\partial \D(c,t)}(v_\nu)
	-M_{\partial \D(c,t)}\bigl(\max\{v_{i,\nu},v_{j,\nu}\}\bigr)
	\le C_{c,s_0,s_1}\frac{\log r_\nu}{r_\nu},
	\]
	since \(\max\{u_i,u_j\}\le u\). Multiplying by \(2\pi t\) and integrating in \(t\) yields
	\begin{equation*}
	0\le
	\int_{A(c;s_0,s_1)}
	\left[v_\nu-\max\{v_{i,\nu},v_{j,\nu}\}\right]\mathrm{d}A
	\le C'_{c,s_0,s_1}\frac{\log r_\nu}{r_\nu}.
	\end{equation*}
	Passing to the limit via \(L^1_{\text{loc}}\)-convergence, it gives that
	\[
	\int_{A(c;s_0,s_1)}
	\left[v-\max\{v_i,v_j\}\right]\mathrm{d}A=0.
	\]
	Such annuli form a topological basis locally, hence the nonnegative difference \(v-\max\{v_i,v_j\}\) vanishes almost everywhere on \(2\D\). The two subharmonic functions \(v\) and \(\max\{v_i,v_j\}\) therefore have the same upper semicontinuous representative, proving \eqref{eq:pair-limit}. Equivalently, it follows that at every point in $2\D$, at least $q-1$ of the functions $v_i$ agree with $v$. 
	
	For $1\le i\le q$, we put \(w_i=v-v_i\ge0\) which is $\delta$-subharmonic and lower semi-continuous, and put \(E_i=\{v=v_i\}\). By Lemma~\ref{lem:grishin}, \(\mathbf1_{E_i}\Delta w_i\ge0\) (i.e. \(\mathbf1_{E_i}\Delta v_i\le \mathbf1_{E_i}\Delta v\)). The restriction of $\Delta w_i$ to $E_i$ is $L^\infty$, since $\Delta v_i$ is positive and $\Delta v$ is $L^\infty$ on $E_i$. Hence Lemma~\ref{lem:dcd-zero-set} gives that
	\begin{equation*}
	\mathbf1_{E_i}\Delta v=\mathbf1_{E_i}\Delta v_i.
	\end{equation*}
	Since $\sum_i\mathbf1_{E_i}\ge q-1$, we conclude that
	\begin{align}
	(q-1)\Delta v
	&\le \sum_{i=1}^q\mathbf1_{E_i}\Delta v
	=\sum_{i=1}^q\mathbf1_{E_i}\Delta v_i
	\le \sum_{i=1}^q\Delta v_i.
	\label{eq:limit-measure-ineq}
	\end{align}
	Let $0<\varepsilon<1/2$ and choose $\chi\in C_c^\infty(\D)$ with $0\le\chi\le1$ and $\chi\equiv1$ on $\D_{1-\varepsilon}$. Testing \eqref{eq:limit-measure-ineq} against $\chi$ and using the weak convergence of \(\{v_\nu\}_{\nu}\) and \(\{v_{i,\nu}\}_{\nu}\), it gives that
	\begin{equation}\label{eq:finite-cutoff}
		(q-1)\int_{\D_{r_\nu}}\chi\left(\frac{z}{r_\nu}\right) \Delta u 
		\le \sum_{i=1}^q \int_{\D_{r_\nu}}\chi\left(\frac{z}{r_\nu}\right) f^*[L_i] 
		+ o(r_\nu^2) + \sum_{i=1}^q \int_{\C}\chi\left(\frac{z}{r_\nu}\right) \mathrm{d}\beta_{i,f}.
	\end{equation}
	Since \(\chi\) is compactly supported, Stokes' theorem gives that
	\[
	\int_{\C}\chi\left(\frac{z}{r_\nu}\right) \mathrm{d}\beta_{i,f}
	= -\int_{\C} \mathrm{d}\left[\chi\left(\frac{z}{r_\nu}\right)\right]\wedge \beta_{i,f}
	= O(r_\nu) = o(r_\nu^2),
	\]
	Using \(\int_{\D_{r_\nu}}\chi(\frac{z}{r_\nu}) \Delta u \ge a_f((1-\varepsilon)r_\nu)\) and \(\int_{\D_{r_\nu}}\chi(\frac{z}{r_\nu}) f^*[L_i] \le n_f(r_\nu, L_i)\), we obtain from \eqref{eq:finite-cutoff} that
	\[
	(q-1)a_f((1-\varepsilon)r_\nu)
	\le \sum_{i=1}^q n_f(r_\nu,L_i) + o(r_\nu^2).
	\]
	Finally, the Brody condition gives a uniform density bound for \(f^*\omega_{\FS}\), hence
	\[
	a_f(r_\nu)-a_f((1-\varepsilon)r_\nu)
	\le C\varepsilon r_\nu^2.
	\]
	Combining with the inequality above and dividing by \(r_\nu^2\), it proves \eqref{eq:limsup} when \(\nu\to\infty\) and \(\epsilon\to 0\).
\end{proof}
\begin{rmk}
	The proof of \eqref{eq:pair-limit} differs from \cite[pp. 1599]{Dacosta-Duval2013}, since our functions \(v_{i,\nu}\) are {not} subharmonic. The subharmonic condition allows one to directly pass to the limit in the circular means. Here, we can only rely on \(L^1_{\mathrm{loc}}\)-convergence. Hence, we use the annulus integral argument to justify the limit.
\end{rmk}

\section{Proof of Theorem \ref{thm:gt-linear}}\label{sec:gt-proof}

\subsection*{Step 1: A preliminary logarithmic growth for \texorpdfstring{$T_f(r)$}{T\_f(r)}}
\label{subsec:coarse-growth}
We begin the proof of Theorem~\ref{thm:gt-linear} with the coarse estimate $O(r\log r)$, as a direct application of Theorem \ref{thm:main}. 

Let $f:\mathbb{C}\to (\mathbb{CP}^2,J)$ be a Brody $J$-holomorphic curve, omitting two distinct $J$-lines $L_1,L_2$. Define a smooth function
\[
v\coloneq u_1-u_2=(\log \|\sigma_1\|-\log \|\sigma_2\|)\circ f.
\]
Since both lines are omitted, \eqref{eq:PL-line-source} yields that  \(\mathrm{d}\mathrm{d}^cv=\mathrm{d}(\beta_{1,f}-\beta_{2,f})\). Thus, up to the fixed normalization constant between $\mathrm{d}\mathrm{d}^c$ and $\Delta$, the function $v$ satisfies that
\begin{equation}\label{eq:dyadic-pde}
		-\Delta v=\operatorname{Div}F,
		\qquad \|F\|_{L^\infty}\le C\|\mathrm{d}f\|\le A.
\end{equation}
Again  for \(i=1,2\), since \(N_f(r,L_i)=0\), \eqref{eq:first-main-line} yields \(M_r(u_i)=O(r)\) where the constant depends on \(L_i\). Hence we have
\begin{equation}\label{eq:dyadic-mean-growth}
M_r(v)=O(r).
\end{equation}
On the other hand, combining together  Proposition~\ref{prop:two-potential}, \eqref{eq:jensen-u} and the equality \(\max\{u_1,u_2\}=\frac12\bigl(u_1+u_2+|v|\bigr)\), we have
\begin{align*}
M_r(|v|)
&= 2M_r(\max\{u_1,u_2\})-M_r(u_1)-M_r(u_2)\\
&= 2\,M_r(u)+O\bigl(r\log(2+r)\bigr)\\
&= 2\,T_f(r)+O\bigl(r\log(2+r)\bigr).
\end{align*}
Now combining with Theorem \ref{eq:main-smt}, it  derives that
	\begin{equation}\label{eq:dyadic-L1-growth}
		M_r(|v|)= o(r^2).
	\end{equation}
\begin{lem}\label{lem:dyadic-harmonic}
	Let $v\in \mathcal{C}^1_{\mathrm{loc}}(\C)$ be a real-valued function, satisfying \eqref{eq:dyadic-pde}-\eqref{eq:dyadic-L1-growth}. Then 
	\begin{equation*}
		v(z)=O\bigl(|z|\log(2+|z|)\bigr),
	\end{equation*}
wherethe implicit constant depends only on the constants in \eqref{eq:dyadic-pde}-\eqref{eq:dyadic-L1-growth}.
\end{lem}

%

Taking Lemma \ref{lem:dyadic-harmonic} for granted, we can strengthen the estimate that
\[
M_r(|v|)=O\bigl(r\log(2+r)\bigr).
\]
Equivalently , we have
\begin{equation}\label{eq:coarse-characteristic}
T_f(r)=O\bigl(r\log(2+r)\bigr).
\end{equation}
Since $a_f(t)$ is nondecreasing and
\[
(\log2)\,a_f(r)
\le \int_r^{2r}a_f(t)\frac{dt}{t}
=T_f(2r)-T_f(r)
\le T_f(2r),
\]
it follows that
\begin{equation}\label{eq:coarse-area}
a_f(r)=O\bigl(r\log(2+r)\bigr).
\end{equation}
\begin{proof}[Proof of Lemma \ref{lem:dyadic-harmonic}]
	For $R>1$, let $q_R$ be the unique weak solution of the Dirichlet problem
	\[
	-\Delta q_R=\operatorname{Div}F\quad\text{in }\D_R,
	\qquad q_R=0\quad\text{on }\partial\D_R.
	\]
	Using the Green function $G_R$ of the disk $\D_R$, which satisfies
	\[-\Delta_\zeta G_R(z,\zeta)=\delta_z(\zeta),\quad G_R(z,\zeta)=0\quad\text{on }\partial\D_R,
	\]
then integrating by parts, we obtain
	\begin{align*}
	q_R(z) &= \int_{\D_R} \delta_z(\zeta) q_R(\zeta) \, \mathrm{d}A(\zeta) 
	= -\int_{\D_R} (\Delta_\zeta G_R(z,\zeta)) q_R(\zeta) \, \mathrm{d}A(\zeta) \\
	&= -\int_{\D_R} G_R(z,\zeta) \Delta_\zeta q_R(\zeta) \, \mathrm{d}A(\zeta) 
	= -\int_{\D_R} G_R(z,\zeta) \operatorname{Div}F(\zeta) \, \mathrm{d}A(\zeta) \\
	&= \int_{\D_R} \nabla_\zeta G_R(z,\zeta) \cdot F(\zeta) \, \mathrm{d}A(\zeta) - \int_{\partial\D_R} G_R(z,\zeta) (F(\zeta)\cdot \nu) \, \mathrm{d}s(\zeta) \\
	&= \int_{\D_R} \nabla_\zeta G_R(z,\zeta) \cdot F(\zeta) \, \mathrm{d}A(\zeta).
	\end{align*}
	Similar to \eqref{eq:green-kernel-estimate}, the Green-kernel estimate gives
	\begin{equation}\label{eq:qR-linear}
		\|q_R\|_{L^\infty(\D_R)}\le O(R).
	\end{equation}
	Put $h_R=v-q_R$ on \( \D_R\). Then $\Delta h_R = 0$, which is harmonic. Since \(q_R=0\) on \(\partial\D_R\), we have $h_R=v$ on $\partial\D_R$.
	
	Take dyadic radii $R_n=2^nR_0$, with \(n\in \mathbb{N}\) and $R_0$ large enough. On each \(\D_{R_n}\), we write the harmonic expansion as
	\[
	h_{R_n}(re^{i\theta})=a_{0,n}+\Re\sum_{k\ge1}a_{k,n}r^ke^{ik\theta}.
	\]
	By the Fourier coefficient formula, for \(k\ge1\), we have
	\[
	a_{k,n}R_n^k =\frac{1}{2\pi}\int_0^{2\pi} h_{R_n}(R_ne^{i\theta})e^{-ik\theta}= \frac{1}{2\pi}\int_0^{2\pi} v(R_ne^{i\theta})e^{-ik\theta}\,\mathrm{d}\theta.
	\]
	Applying  \eqref{eq:dyadic-L1-growth}, it gives that for \(k\ge1\),
	\begin{equation}\label{eq:coeff-rough}
		|a_{k,n}|\le \varepsilon_n R_n^{2-k},
	\end{equation}
	where \(\varepsilon_n \to 0\) as \(n \to \infty\). For $k=0$, \eqref{eq:dyadic-mean-growth} and the mean formula gives that
	\begin{equation}\label{eq:a0-bound}
		|a_{0,n}|=\left|\frac1{2\pi}\int_0^{2\pi}h_{R_n}(R_ne^{i\theta})\,d\theta\right|=\left|\frac1{2\pi}\int_0^{2\pi}v(R_ne^{i\theta})\,d\theta\right|\le C R_n.
	\end{equation}
	
On $\D_{R_n}$, the difference of two consecutive harmonic approximations satisfies
\[\left|h_{R_{n+1}}-h_{R_n}\right|_{L^\infty(\D_{R_n})}=\left|q_{R_n}-q_{R_{n+1}}\right|_{L^\infty(\D_{R_n})}\le C R_n.\] 
Thus by Cauchy's estimates for harmonic functions, we have, for \(k\ge1\),
	\begin{equation}\label{eq:coeff-difference}
		|a_{k,n+1}-a_{k,n}|\le C R_n^{1-k}.
	\end{equation}
For every \(k\ge2\), \eqref{eq:coeff-rough} implies that \(a_{k,n}\to0\) as \(n\to\infty\). Summing \eqref{eq:coeff-difference} from \(n\) to infinity, it gives the sharper bound for all \(k\ge2\) that
	\begin{equation}\label{eq:high-coeff}
		|a_{k,n}|\le\sum_{m=n}^\infty|a_{k,m+1}-a_{k,m}|\le \sum_{m=n}^\infty C R_m^{1-k}\le C R_n^{1-k}.
	\end{equation}
For \(k=1\), it gives that
	\begin{equation}\label{eq:linear-coeff}
		|a_{1,n}|\le |a_{1,0}|+\sum_{m=0}^{n-1}C R_m \le C(1+n)\le C\,\log(2+R_n).
	\end{equation}
	
For an arbitrary \(z\), there exists \(R_n\) such that \(|z|\le \frac{R_n}{4}\). Then we combine \eqref{eq:a0-bound}-\eqref{eq:linear-coeff} to estimate \(h_{R_n}(z)\) as
	\begin{align}\label{eq:hRn-log}
	|h_{R_n}(z)|
	&\le |a_{0,n}|+|a_{1,n}||z|+\sum\limits_{k\ge2}|a_{k,n}||z|^k\nonumber\\
	&\le C R_n+C\log(2+R_n)\cdot \frac{R_n}{4}+\sum\limits_{k\ge2} C R_n^{1-k} \left(\frac{R_n}{4}\right)^k\nonumber\\
	&\le C R_n\log(2+R_n).
\end{align}
Combining this with \eqref{eq:qR-linear}, we then obtain \(|v(z)|=|h_{R_n}(z)+q_{R_n}(z)|\le C R_n\log(2+R_n)\). For any large \(|z|\), choose \(n\) such that \(4|z|\le R_n<8|z|\), then the bound becomes \(O(|z|\log(2+|z|))\). This proves the desired result.
\end{proof}

\begin{rmk}
	In Lemma \ref{lem:dyadic-harmonic}, if we replace \eqref{eq:dyadic-L1-growth} by 
	\begin{equation}\label{eq:dyd-L1-Or2}
	M_r(|v|)= O(r^2),
	\end{equation}
then the conclusion of Lemma~\ref{lem:dyadic-harmonic} weakens to the existence of a real homogeneous harmonic quadratic polynomial $Q$ such that
\begin{equation}
\label{eq:dyadic-conclusion-2}
v(z)=Q(z)+O\bigl(|z|\log(2+|z|)\bigr).
\end{equation}
Indeed, under assumption \eqref{eq:dyd-L1-Or2}, when \(k=2\), the Fourier coefficient estimate \eqref{eq:coeff-rough} yields \(|a_{2,n}|\le C\). Then the estimate \eqref{eq:coeff-difference} implies that \(\{a_{2,n}\}_n\) is a Cauchy sequence, converging to some limit \(a_2\) with \(|a_{2,n}-a_2|\le \frac{C}{R_n}\). Let \(Q(z)=\Re(a_2z^2)\). Similarly, we have
\[
|h_{R_n}(z)-Q(z)|\le C R_n\log(2+R_n),
\]
therefore \eqref{eq:dyadic-conclusion-2} holds. It is the quantitative substitute, in the almost complex setting, for the fact that an order-two holomorphic curve is the exponential of a polynomial of degree at most two.
\end{rmk}
\subsection*{Step 2: Linear growth of \texorpdfstring{$M_r(|v|)$}{M\_r(|v|)}}\label{sec:weighted-splitting}
In the second step, we shall apply the weighted \(L^2\)-estimate in functional analysis to globalize the proof of Lemma \ref{lem:dyadic-harmonic} and sharpen the estimate \eqref{eq:qR-linear}. Therefore, we can obtain a global and linear version of \eqref{eq:hRn-log}, which implies a linear growth of $M_r(|v|)$. 

Recall that the smooth function \(v=u_1-u_2\) satisfies \eqref{eq:dyadic-pde}. Fix $0<\varepsilon<1$. Set $\beta=1+\varepsilon\in(1,2)$ and $\alpha=\frac\beta2\in\left(\frac12,1\right)$. Define the weight function
\begin{equation*}
\mathfrak w(z)=(1+|z|)^{-\beta},
\end{equation*}
which belongs to the Muckenhoupt class $A_2(\mathbb R^2)$. That is, for any \(a\in \mathbb{C}\) and any \(r>0\), it satisfies
\begin{equation}\label{eq:Muckenhoupt-2}
	\left( \frac{1}{|\D(a,r)|}\int_{\D(a,r)} \mathfrak w \,\mathrm{d}A \right) \left( \frac{1}{|\D(a,r)|}\int_{\D(a,r)} \mathfrak w^{-1} \,\mathrm{d}A \right) \le C.
\end{equation} 
Indeed, when \(r\le\frac{1+|a|}{2}\), we have \(\mathfrak w(z) \asymp C (1+|a|)^{-\beta}\) for all \(z \in \D(a,r)\). It is comparable to a constant on the disk, hence \eqref{eq:Muckenhoupt-2} holds. When \(r > \frac{1+|a|}{2}\), there exists a constant \(C\) such that \(\D(a,r)\subset\D(0,Cr)\). Since \(\beta < 2\), we have
	\[
	\int_{\D(a,r)} \mathfrak w \,\mathrm{d}A \le \int_{\D(0,Cr)} (1+|z|)^{-\beta} \,\mathrm{d}A \le C \int_0^{Cr} (1+s)^{-\beta} s \,\mathrm{d}s \le C r^{2-\beta},
	\]
and
	\[
	\int_{\D(a,r)} \mathfrak w^{-1} \,\mathrm{d}A \le C \int_0^{Cr} (1+s)^\beta s \,\mathrm{d}s \le C r^{2+\beta}.
	\]
Since \(|\D(a,r)|=O(  r^2)\), we obtain \eqref{eq:Muckenhoupt-2}. And since the constant \(C\) is uniform in both cases, we conclude that \([\mathfrak w]_{A_2}<\infty\).

\begin{lem}\label{lem:weighted-splitting}
	There exist a function $q\in W^{1,2}_{\rm loc}(\C)$ and a harmonic function $h\in \mathcal{H}(\C)$ such that \(v\) has the decomposition as 
	\begin{equation*}
	v=h+q,
	\qquad -\Delta q=\operatorname{Div}F.
	\end{equation*}
	And after adding a constant to $q$, we have
	\begin{equation*}
	\frac1{2\pi}\int_0^{2\pi}|q(Re^{i\theta})|\,\mathrm{d}\theta
	\le C\bigl(1+R^\alpha\bigr),\qquad R\ge1.
	\end{equation*}
\end{lem}

\begin{proof}
Firstly we have \(F\in L^2(\mathfrak w)\), since
\begin{align}
	\int_{\C}|F|^2\mathfrak w\,\mathrm{d}A
	&= \int_{\D}|F|^2\mathfrak w\,\mathrm{d}A+\sum_{k\ge1}\int_{\D_{2^{k}}\setminus\D_{2^{k-1}}}|F|^2\mathfrak w\,\mathrm{d}A	\notag\\
	&\le \int_{\D}|F|^2\mathfrak w\,\mathrm{d}A
	+ C\sum_{k\ge1}2^{-\beta k}\int_{\D_{2^{k}}\setminus\D_{2^{k-1}}}|\mathrm{d}f|^2\,\mathrm{d}A	&&\text{(by $|F|\le C|\mathrm{d}f|$ in \eqref{eq:dyadic-pde})}	\notag\\
	&\le C + C\sum_{k\ge0}2^{-\beta k}\int_{\D_{2^{k+1}}}|\mathrm{d}f|^2\,\mathrm{d}A	\notag\\
	&\le C + C\sum_{k\ge0}2^{-\beta k}\cdot\bigl(2^{k+1}\log(2+2^{k+1})\bigr)
	&&\text{(by \eqref{eq:coarse-area})}
	\notag\\
	&\le C\sum_{k\ge0}(k+1)2^{(1-\beta)k}<\infty.	&&\text{(since $\beta>1$)}
	\label{eq:F-weighted}
\end{align}

Let $\mathcal P(F)=\nabla\Delta^{-1}\operatorname{Div}(F)$ be the gradient projection, which is a Calder\'on--Zygmund operator. Its matrix entries are compositions of Riesz transforms. Since the operator $\mathcal P$ is initially defined via global Fourier multipliers, we regularize the domain by choosing a sequence of cutoff functions $\{\chi_N\}_{N\ge 1}\in C_c^\infty(\mathbb{C})$ with $0\le\chi_N\le1$, $\chi_N=1$ on $\D_N$ and $\chi_N=0$ outside $\D_{2N}$. Setting $F_N=\chi_NF$, it follows from \eqref{eq:F-weighted} and the dominated convergence theorem that
	\[
	F_N\longrightarrow F\quad\text{in }L^2(\mathfrak w).
	\]
	Since \(\mathfrak w\in A_2(\mathbb R^2)\), the weighted \(L^2\)-theorem for Calder\'on--Zygmund operators \cite{Hytonen2012} states that there exists a constant \(C>0\) depending only on \(\mathfrak w\) such that, for any vector field \(U \in L^2(\mathfrak w)\)
		\[
		\|\mathcal P(U)\|_{L^2(\mathfrak w)} \le C \|U\|_{L^2(\mathfrak w)}.
		\]
	Set $G_N=\mathcal P(F_N)$. Consequently, there exists a constant \(C>0\) such that
		\[
		\|G_N-G_M\|_{L^2(\mathfrak w)}
		= \|\mathcal P(F_N - F_M)\|_{L^2(\mathfrak w)}
		\le C\|F_N-F_M\|_{L^2(\mathfrak w)}.
		\]
	Hence \(\{G_N\}_N\) is also a Cauchy sequence in \(L^2(\mathfrak w)\). By the completeness of \(L^2(\mathfrak w)\), there exists a vector field \(G \in L^2(\mathfrak w)\) such that \(G_N \to G\) strongly in \(L^2(\mathfrak w)\), with
	\begin{equation}\label{eq:weighted-Hodge}
	\|G\|_{L^2(\mathfrak w)}\le C\|F\|_{L^2(\mathfrak w)}.
	\end{equation}
	For every $N$, Fourier multiplier identities give that, in the sense of distributions,
	\[
	\operatorname{curl}G_N=0,
	\qquad
	\operatorname{Div}G_N=\operatorname{Div}F_N.
	\]
	Since $\mathfrak w$ is bounded below by a positive constant on each compact set, strong convergence in $L^2(\mathfrak w)$ implies local $L^2$-convergence, hence convergence in $\mathcal{D}'(\mathbb{C})$. Testing the above identities against an arbitrary $\varphi\in C_c^\infty(\mathbb{C})$ and passing to the limit, using that $\nabla\varphi$ is compactly supported and that $G_N\to G$, $F_N\to F$ in $L^2_{\mathrm{loc}}$, we obtain
\[
\operatorname{curl}G=0,
\qquad
\operatorname{Div}G=\operatorname{Div}F
\quad\text{in }\mathcal{D}'(\mathbb{C}).
\] 
As \(\mathbb{C}\) is simply connected and $G$ is curl-free in the distributional sense, the weak Poincar\'e lemma yields a potential $q\in W^{1,2}_{\rm loc}(\mathbb{C})$, unique up to an additive constant, with $G=\nabla(-q)$. Then
\[
-\Delta q=\operatorname{Div}G=\operatorname{Div}F
\quad\text{in }\mathcal{D}'(\mathbb{C}).
\]
And since also $-\Delta v=\operatorname{Div}F$ , we get $\Delta h=0$ in $\mathcal{D}'(\mathbb{C})$ for $h=v-q$. By Weyl's lemma, $h\in C^\infty(\mathbb{C})$ and $\Delta h=0$ in the classical sense, i.e. $h$ is harmonic on $\mathbb{C}$ \cite{Hormander03}.
	
	On $\D_{2R}$, one has $\mathfrak w(z)\ge C R^{-\beta}$ for $R\ge1$.  Hence \eqref{eq:weighted-Hodge} gives that
	\[\int_{\D_{2R}}|\nabla q|^2\,\mathrm{d}A \cdot R^{-\beta}
	\le C \int_{\D_{2R}}|\nabla q|^2 \mathfrak{w}\,\mathrm{d}A
	\le C \int_{\C}|F|^2 \mathfrak{w}\,\mathrm{d}A\le C.\] That is,
	\begin{equation}\label{eq:q-energy-growth}
		\int_{\D_{2R}}|\nabla q|^2\,\mathrm{d}A
		\le C R^\beta.
	\end{equation}
	On the annulus \(A_R = \{R/2 < |z| < 2R\}\), put \(c_R = \frac{1}{|A_R|} \int_{A_R} q \,\mathrm{d}A\).
	Since the diameter of \(A_R\) is comparable to \(R\), Poincar\'e's inequality implies that
	\begin{equation}\label{eq:Poin}
		\int_{A_R}|q-c_R|^2\,\mathrm{d}A \le C R^2 \int_{A_R}|\nabla q|^2\,\mathrm{d}A.
	\end{equation}
	Combining this with \eqref{eq:q-energy-growth} and dividing by \(|A_R| \asymp R^2\), we obtain
	\begin{equation}\label{eq:campanato-annulus}
		\left(\frac1{|A_R|}\int_{A_R}|q-c_R|^2\,\mathrm{d}A\right)^{1/2}
		\le C R^{\beta/2} = C R^\alpha.
	\end{equation}
	The standard \(W^{1,2}\) trace theorem on the annulus \(A=\{1/2 < |z| < 2\}\) states that, for \(f\in W^{1,2}(A)\),
		\[
		\int_{\partial\D} |f|^2\,\mathrm{d}s
		\le C\left[ \int_{A} |f|^2\,\mathrm{d}A+ \int_{A}|\nabla f|^2\,\mathrm{d}A \right].
		\]
	We apply this to the scaled function \(f(w)=q(Rw)-c_R\) on the annulus \(A_R = \{R/2 < |z| < 2R\}\). Then we obtain
		\begin{equation*}
		\frac1{2\pi R}\int_{\partial\D_R}|q-c_R|^2\,\mathrm{d}s
		\le C\left[
		R^{-2}\int_{A_R}|q-c_R|^2\,\mathrm{d}A
		+\int_{A_R}|\nabla q|^2\,\mathrm{d}A
		\right].
		\end{equation*}
	Combining this \eqref{eq:Poin} and \eqref{eq:q-energy-growth}, we obtain
	\begin{align}
	\left(\frac1{2\pi R}\int_{\partial\D_R}|q-c_R|^2\,\mathrm{d}s\right)^{1/2}\le\left( \int_{A_R}|\nabla q|^2\,\mathrm{d}A \right)^{1/2} \le C R^\alpha. \label{eq:q-trace}
	\end{align}
	The annuli \(A_R\) and \(A_{2R}\) overlap on a set of area asymptotic to \(R^2\). Also by the triangle inequality and \eqref{eq:campanato-annulus}, we obtain \(|c_{2R}-c_R| \le CR^\alpha\). For arbitrary \(R\ge1\), we can choose an integer \(m\) such that \(2^m \le R < 2^{m+1}\). Similarly working on the overlap set $A_R\cap A_{2^m}$, it yields that \(|c_R-c_{2^m}|\le CR^\alpha\).
	On the other hand, 
	\[|c_{2^m}-c_1|\le\sum\limits_{k=0}^{m-1}|c_{2^{k+1}}-c_{2^k}|\le C2^{m\alpha}\le CR^\alpha.\]
	It implies that \(|c_R-c_1|\le CR^\alpha\) for every \(R\ge1\). Absorbing \(c_1\) into the additive normalization of \(q\), \eqref{eq:q-trace} and Cauchy--Schwarz inequality imply that
	\begin{align*}
		M_R(|q|) 
		&\le M_R(|q-c_R|) + |c_R| \\
		&\le C \left(\frac1{2\pi R}\int_{\partial\D_R}|q-c_R|^2\,\mathrm{d}s\right)^{1/2} + CR^\alpha \\
		&\le C(1+R^\alpha).
	\end{align*}
\end{proof}

\begin{pro}\label{prop:v-linear}
	We have the linear growth that
	\begin{equation}\label{eq:v-linear}
	M_r(|v|)=O(r).
	\end{equation}
\end{pro}

\begin{proof}
	Combining \(M_r(|v|)=O\bigl(r\log(2+r)\bigr)\) with Lemma~\ref{lem:weighted-splitting}, it gives that
	\[
	M_r(|h|)=O\bigl(r\log(2+r)\bigr).
	\]
	Since $h$ is a harmonic function defined on \(\C\), we write $h=\Re H$ with $H=\sum_{m\ge0}a_mz^m$ an entire function. The Fourier coefficient formula gives that
	\[
	|a_m|r^m\le C M_r(|h|)\le Cr\log(2+r).
	\]
	Letting $r\to\infty$, it shows that $a_m=0$ for any $m\ge2$. Thus $h$ is affine and \(M_r(|h|)=O(r)\). Applying Lemma~\ref{lem:weighted-splitting} and since $\alpha<1$, it gives that
\[
M_r(|v|)\le M_r(|h|)+M_r(|q|)=O(r),
\]
as required.
\end{proof}

\subsection*{Step 3: Potential comparison with linear error term}\label{sec:linear-pair}
This step further develops the techniques in Lemma \ref{lem:perturbed} and sharpen Proposition \ref{prop:two-potential} with linear error term.

Let $P=L_1\cap L_2$ and choose an auxiliary $J$-line $L_3$ such that $P\notin L_3$. Recall that $\lambda_3=\log\|\sigma_3\|$ and that \(\widehat\lambda_{12}=\frac12\log\bigl(\|\sigma_1\|_h^2+\|\sigma_2\|_h^2\bigr)\). As \eqref{eq:central-PL}, we have
\begin{equation}\label{eq:duval-P-potential}
\pi_P^*\omega_{\mathbb P^1}
=\omega_{\FS}+\mathrm{d}\mathrm{d}_J^c\psi_{12}+\mathrm{d}a_{12},
\end{equation}
with $\psi_{12}(x)=\log\dist(x,P)+O(1)=\widehat\lambda_{12}(x)+O(1)$ near \(P\), and \(a_{12}\) is bounded on $\mathbb{CP}^2\setminus\{P\}$. Choose a relatively compact coordinate ball $U$ around $P$, whose closure is disjoint from $L_3$.  Since $\omega_{\FS}$ is closed and $U$ is contractible, we can fix a smooth one-form $\vartheta$ on $U$ such that
\begin{equation}\label{eq:local-primitive}
\mathrm{d}\vartheta=\omega_{\FS}\quad\hbox{on }U.
\end{equation}

\begin{lem}\label{lem:target-localization}
	After adding a constant to $\lambda_3$, the target set
	\[
	\mathcal O=\{x\in\mathbb{CP}^2\setminus\{P\}:\lambda_3(x)>\psi_{12}(x)\}
	\]
	is nonempty and satisfies $\overline{\mathcal O}\Subset U$.
\end{lem}
\begin{proof}
The function $\psi_{12}$ is smooth on \(\mathbb{CP}^2\setminus\{P\}\) and has only logarithmic pole at \(P\). Thus, it is bounded from below on the compact set \(\mathbb{CP}^2\setminus U\), whereas $\lambda_3$ is bounded from above there. By subtracting a sufficiently large constant from $\lambda_3$, we obtain \(\lambda_3<\psi_{12}\) on \(\mathbb{CP}^2\setminus U\).

Conversely, since $\psi_{12}\to-\infty$ near \(P\) and $\lambda_3$ is finite near \(P\), we have \(\lambda_3>\psi_{12}\) on a punctured neighborhood of \(P\). Shrinking this neighborhood to lie inside \(U\) establishes the results.
\end{proof}

Set
\[
u_{12} \coloneq u+\psi_{12}\circ f,\qquad
\Omega\coloneq f^{-1}(\mathcal O),\qquad
w \coloneq u_3-u_{12}.
\]
Since the two omitted lines intersect at \(P\), the curve \(f\) omits \(P\); consequently, \(u_{12}\) is a smooth function and \(w>0\) on \(\Omega\). If \(\Omega=\C\), then \(f(\C)\subset\overline{\mathcal O}\Subset U\). By \eqref{eq:local-primitive}, applying Stokes' theorem yields
\begin{equation}
\label{eq:all-in-U}
a_f(r)=\int_{\D_r}f^*\omega_{\FS}
=\int_{\partial\D_r}f^*\vartheta
=O(r),
\end{equation}
where the last estimate follows from the Brody bound. We may therefore assume \(\partial\Omega\ne\varnothing\). By definition, \(\Omega\) does not meet \(f^{-1}(L_3)\). Pulling back \eqref{eq:duval-P-potential} via \(f\) and applying \eqref{eq:PL-line-target}, it therefore yields that on $\Omega$
\[
-\mathrm{d}\mathrm{d}^cw=(\pi_P\circ f)^*\omega_{\mathbb P^1}+\mathrm{d}\eta,
\]
where $\eta=-f^*a_{12}-\beta_{3,f}$ with $\beta_{3,f}$ bounded on $\Omega$. Thus we obtain
\begin{equation}\label{eq:localized-w-equation}
-\Delta w=\mu+\operatorname{Div}F,
\end{equation}
where \(\mu\) is a fixed positive multiple of \((\pi_P\circ f)^*\omega_{\mathbb P^1}\), and
\begin{equation}\label{eq:localized-F-bound}
 \|F\|_{L^\infty(\C)}\le C|\mathrm{d}f|\le C.
\end{equation}
\begin{rmk}
	The inequality \eqref{eq:localized-F-bound} highlights a fundamental difference between the omitted line case and the general setting of Lemma~\ref{lem:perturbed} and Proposition~\ref{prop:two-potential}. In Section~\ref{sec:pair}, where \(f\) might intersect \(P\), the decompositions \eqref{eq:uij-decomp} and \eqref{eq:w-current} show that the divergence term \(F\) originates from \(\tilde{f}^*(\sigma_{12}^*a_{12})\) and \(f^*\beta_{3}\), whose norm never arrives an upper bound from \(|\mathrm{d}f|\) merely. In contrast, the bound \eqref{eq:localized-F-bound} holds precisely because the omitted lines force \(f\) to stay away from the singularity \(P\), so that \(F\) satisfies \(|F| \le C|\mathrm{d}f|\).
\end{rmk}
\begin{lem}\label{lem:linear-energy-maximal}
	Let $D=\D(a,2\rho)\subset\Omega$ satisfy \(2\rho=\dist(a,\partial\Omega)\), and let \(z_0\in\partial{D}\cap\partial\Omega\). Then for every \(0<s\le4\rho\), we have
	\begin{equation}\label{eq:morrey-F}
	\int_{{D}\cap \mathbb{D}(z_0,s)}|F|^2\,\mathrm{d}A\le C s.
	\end{equation}
	The constant \(C\) is independent of \(a,\rho,z_0\) and \(s\).
\end{lem}
\begin{proof}
	Since \(D\subset\Omega\), Lemma~\ref{lem:target-localization} gives \(f(D\cap \mathbb{D}(z_0,s))\subset U\). By \eqref{eq:local-primitive} and Stokes' theorem, we have
	\[
	\int_{D	\cap \mathbb{D}(z_0,s)}f^*\omega_{\FS}=\int_{\partial (D\cap \mathbb{D}(z_0,s))}f^*\vartheta.
	\]
	The intersection \(D\cap \mathbb{D}(z_0,s)\) is the intersection of two disks, and its boundary consists of at most two circular arcs. For \(s\le4\rho\), the total length of these arcs is bounded by \(Cs\). Since \(\vartheta\) is bounded on \(\overline U\) and \(f\) is Brody, we get
	\[
	0\le\int_{D\cap \mathbb{D}(z_0,s)}f^*\omega_{\FS}\le Cs.
	\]
	Finally, by \eqref{eq:localized-F-bound} and the area \(|\mathrm{d}f|^2\,\mathrm{d}A\) is uniformly comparable to \(f^*\omega_{\FS}\), it proves \eqref{eq:morrey-F}.
\end{proof}

The next estimate is the place where the logarithm disappears.

\begin{lem}\label{lem:morrey-normal}
	Let $D=\D(a,2\rho)$ and $z_0\in\partial D$. Let $p$ be a function solving the Dirichlet equation
	\[
	\Delta p=\operatorname{Div}F\quad\hbox{in }D,
	\qquad p=0\quad\hbox{on }\partial D.
	\]
	Assume $\|F\|_\infty\le A$, and assume that for some fixed $r_*>0$, for any \(z\in D\cap \D(z_0,2r_*)\), it holds
	\[
	|F(z)-F(z_0)|\le L|z-z_0|,
	\]
	and
	\[
	\int_{D\cap \D(z_0,s)}|F|^2\,\mathrm{d}A\le C_0s
	\qquad(0<s\le4\rho).
	\]
	Then we have
	\begin{equation*}
	|\partial_\nu p(z_0)|\le C,
	\end{equation*}
	where $\nu$ denotes the inward unit normal to $\partial D$ and the constant  $C$ depends only on $|F|,L,r_*,C_0$.
\end{lem}

\begin{proof}
	For $\rho\le2r_*$, the inequality \eqref{eq:E-bound} in the proof of Lemma~\ref{lem:perturbed} already gives the result. Hence we assume $\rho>2r_*$. Similar to the proof of Lemma~\ref{lem:perturbed}, Green's representation theorem yields
	\[
		p(x) = -\int_{D} \nabla_t G(x,t) \cdot F(t) \,\mathrm{d}t.
	\]
	Let \(P_D(t) = \partial_{\nu} G(z_0,t) = \frac{1}{4\pi\rho}\,\frac{(2\rho)^2-|t-a|^2}{|t-z_0|^2}\) be the Poisson kernel on $D$ with pole at $z_0$. To estimate the normal derivative \(\partial_\nu p(z_0)\), we cannot differentiate the integral directly, because \(\nabla P_D(t)\) has a non-integrable singularity of order \(|t-z_0|^{-2}\) at $z_0$. In Lemma~\ref{lem:perturbed}, this singularity is resolved by replacing \(F(t)\) with \(F(t)-F(z_0)\) and applying the Lipschitz estimate on $D\cap \D(z_0,2r_*)$, at the cost of an \(O(r\log r)\) error term on $D\setminus \D(z_0,2r_*)$ (see \eqref{eq:E-bound}). However, equipped with \eqref{eq:morrey-F}, we adapt a new approach here that eliminates the logarithmic loss on the complement $D\setminus \D(z_0,2r_*)$.
	
	Let \(\chi\) be a smooth cutoff function such that \(\chi\equiv1\) on \(\D(z_0,r_*)\), \(\operatorname{supp}\chi\subset \D(z_0,2r_*)\), and \(|\nabla\chi|\le C/r_*\). In the sense of distributions, integrating by parts yields the normal derivative
	\begin{equation}\label{eq:regularized-normal-pairing}
		\partial_\nu p(z_0) = -\langle\nabla P_D, F\rangle,
	\end{equation}
	where the pairing on the right hand side is defined as follows:
		\begin{align}
		\langle\nabla P_D,F\rangle
		&\coloneq\int_D\chi\,\nabla P_D\cdot(F-F(z_0))
		+\int_D(1-\chi)\,\nabla P_D\cdot F \notag\\
		&\qquad
		+F(z_0)\cdot\langle\nabla P_D,\chi\rangle.
		\label{eq:regularized-pairing-decomp}
		\end{align}
The first term is bounded by the Lipschitz estimate \(|F(t)-F(z_0)|\le L|t-z_0|\) and \(\operatorname{supp}\chi\subset \D(z_0,2r_*)\). More precisely, \(\left|\int_D\chi\,\nabla P_D\cdot(F-F_0)\right| \le CLr_*\). The third term can be identified as 
\[\left|F(z_0)\cdot\langle\nabla P_D,\chi\rangle\right|=\left|F(z_0)\cdot(-\int_D P_D\nabla\chi)\right|\le C|F|,\]
since on \(\operatorname{supp}(\nabla\chi)\subset\D(z_0,2r_*)\setminus\D(z_0,r_*)\), we have \(P_D\le C/r_*\).

For the second term, \(1-\chi\equiv0\) on \(\D(z_0,r_*)\) and the integration domain $D\setminus \D(z_0,r_*)$ is separated from \(z_0\). Decompose its support into annuli
	\[
	A_j=D\cap\{s_j\le|t-z_0|<2s_j\},
	\]
where $s_j=2^jr_*$ with $j\in \mathbb{Z}_{\ge 0}$ up to the diameter of \(D\). For every annulus except possibly the last one, we have \(2s_j\le4\rho\). And Lemma \ref{lem:linear-energy-maximal} together with Cauchy-Schwarz gives that
	\[
	\int_{A_j}|F|\,\mathrm{d}A
	\le |A_j|^{1/2}
	\left(\int_{D\cap \D(z_0,2s_j)}|F|^2\,\mathrm{d}A\right)^{1/2}
	\le Cs_j^{3/2}.
	\]
	For the last annulus, where \(2s_j\) may exceed \(4\rho\), we use the whole-disk estimate that
	\[
	\int_D|F|^2\,\mathrm{d}A\le C\rho\le Cs_j.
	\]
Thus \(\int_{A_j}|F|\,\mathrm{d}A\le Cs_j^{3/2}\) for every \(j\). Since \(|\nabla P_D|\le Cs_j^{-2}\) on \(A_j\), we get
	\[
	\int_{A_j}|\nabla P_D|\,|F|\,\mathrm{d}A
	\le Cs_j^{-2}s_j^{3/2}=Cs_j^{-1/2}.
	\]
	The series \(\sum_j s_j^{-1/2}\) converges uniformly, providing a uniform bound for \(\int_D(1-\chi)\,\nabla P_D\cdot F\).

	Combining these estimates with \eqref{eq:regularized-normal-pairing} and \eqref{eq:regularized-pairing-decomp}, we finish the proof.
\end{proof}
Now we are capable to withdraw the logarithm error in pair comparison.
\begin{pro}\label{prop:linear-pair}
	If \(f\) omits \(L_1\) and \(L_2\), then
	\begin{equation}\label{eq:linear-pair}
	M_r(u)=M_r\bigl(\max\{u_1,u_2\}\bigr)+O(r).
	\end{equation}
\end{pro}
\begin{proof}
	If $\Omega=\C$, the stronger conclusion $T_f(r)=O(r)$ follows from \eqref{eq:all-in-U}; hence we assume $\partial\Omega\ne\varnothing$. We apply the technique in the proof of Proposition~\ref{prop:two-potential} to \eqref{eq:localized-w-equation}--\eqref{eq:localized-F-bound}. For any $a\in\Gamma_r\cap\Omega$, set
	\[
	2\rho(a)=\dist(a,\partial\Omega),\quad
	2D_a=\D(a,2\rho(a)),\quad D_a=\D(a,\rho(a)).
		\]
	Choose \(r\) appropriately large and translate \(\Gamma_r\) if necessary such that \(\rho(a)\le r\). Choose \(z_0\in\partial (2D_a)\cap\partial\Omega\). 
		
	The contact set
		\[
		K_P:=\{x\in\mathbb{CP}^2\setminus\{P\}:\lambda_3(x)=\psi_{12}(x)\}
		\]
	is contained in a compact subset of \(\mathbb{CP}^2\setminus(L_3\cup\{P\})\). Hence \(\mathrm{d}\lambda_3\), \(\mathrm{d}\psi_{12}\), \(a_{12}\), and \(\beta_3\) are smooth with uniform bounds on a fixed neighborhood of \(K_P\). Exactly as in \eqref{eq:w-scalar}--\eqref{eq:F-contact-lip}, the Brody condition and interior elliptic estimates give a fixed \(r_*>0\), a uniform bound on contact gradient, and a uniform Lipschitz bound for the coefficients of \(F\) near every contact point \(z_0\). 
		
We now apply Lemma~\ref{lem:linear-energy-maximal} and Lemma~\ref{lem:morrey-normal} as follows. Let \(p_a\) be the solution of the Dirichlet problem
\[
	\Delta p_a=\operatorname{Div}F\quad\text{on }2D_a,
	\qquad p_a|_{\partial(2D_a)}=0.
	\]
Then Lemma~\ref{lem:morrey-normal} gives \(|\partial_\nu p_a(z_0)|\le C\), while the standard Green-kernel estimate gives \(\|p_a\|_{L^\infty(2D_a)}\le C\rho(a)\). On \(2D_a\), we have
	\[
	-\Delta w=\widetilde\mu+\operatorname{Div}F,
	\qquad \widetilde\mu\ge0.
	\]
And we set \(Q_a:=w+p_a\) such that
	\[
	-\Delta Q_a=\widetilde\mu\ge0,
	\qquad Q_a|_{\partial(2D_a)}=w|_{\partial(2D_a)}\ge0,
	\qquad Q_a(z_0)=0.
	\]
Hence \(Q_a\ge0\) on \(2D_a\) by the minimum principle. Moreover the contact gradient bound and Lemma~\ref{lem:morrey-normal} give that
	\[
	D_\nu^+Q_a(z_0)
	\le |\nabla w(z_0)|+|\partial_\nu p_a(z_0)|\le C.
	\]
We apply \cite[Lemma 2]{Dacosta-Duval2013} to obtain that
	\[
	M_{\Gamma_r\cap D_a}(Q_a)\le C\rho(a).
	\]
Hence for \(w=Q_a-p_a\) and \(\|p_a\|_{L^\infty(2D_a)}\le C\rho(a)\), we have \(M_{\Gamma_r\cap D_a}(w)\le C\rho(a).\)
	
Again, we apply the Besicovitch covering theorem to obtain a subcover $\{D_{a_k}\}_k$ satisfying \(\sum\limits_k\rho_k\le Cr\). Summing over \(k\), we obtain \(M_r(w_+)=O(r)\). Recall that \(u=\max\{u_1,u_2\}+w_++O(1)\),	then we finish the proof.
\end{proof}

\subsection*{Step 4: End of the proof}
	Since $L_1,L_2$ are omitted, the FMT \eqref{eq:first-main-line} gives that
	\[
	M_r(u_1)=O(r),\qquad M_r(u_2)=O(r).
	\]
 Combining \(\max\{u_1,u_2\}=\frac12\bigl(u_1+u_2+|v|\bigr)\) and Proposition~\ref{prop:v-linear}, it yields
	\[
	M_r\bigl(\max\{u_1,u_2\}\bigr)=O(r).
	\]
	Therefore \eqref{eq:linear-pair} gives that
	\[
	T_f(r)=M_r(u)=O(r).
	\]
	Finally, by the same arguments as \eqref{eq:coarse-characteristic}-\eqref{eq:coarse-area}, we also get \eqref{eq:gt-linear-area}.
\section*{AI disclosure} The first draft of this paper, which was completely written by the authors, consists of the full proof of Theorem \ref{thm:main} and a weaker form of Theorem~\ref{thm:gt-linear}, in which the growth rate $O(r\log r)$ was obtained. To attain the linear growth rate $O(r)$, AI (GPT Sol 6.1) was used to assist with technical estimates (Lemma~\ref{lem:weighted-splitting} and Lemma~\ref{lem:morrey-normal}). AI tools were also used for searching the literature and polishing the writing. The authors have independently verified the mathematical arguments and  take full responsibility for the accuracy and content of the paper.
\section*{Acknowledgements}
 D.T. Huynh is supported by the Vietnam National Foundation for Science and Technology Development (NAFOSTED) under the grant number 101.04-2025.19. Yunling Chen thanks to the Academy of Mathematics and Systems Science in Beijing for the great research environment.
\begin{center}
	\bibliographystyle{alpha}
	\bibliography{article, references}
\end{center}

\end{document}